\documentclass[10pt]{amsart}
\usepackage[english]{babel}
\usepackage[letterpaper,top=3.25cm,bottom=2cm,left=3cm,right=3cm,marginparwidth=1.75cm]{geometry}

\usepackage{graphicx} 
\usepackage{amsthm, amsmath, amssymb} 
\usepackage{physics}
\usepackage{microtype}
\usepackage{setspace}
\usepackage{subfiles}
\usepackage{color,soul}
\usepackage{tikz-cd} 
\usepackage{mathtools}
\usepackage{enumerate}

\theoremstyle{plain}
   \newtheorem{theorem}{Theorem}[section]
   \newtheorem{prop}[theorem]{Proposition}
   \newtheorem{lemma}[theorem]{Lemma}
   \newtheorem{coro}[theorem]{Corollary}
   \newtheorem{definition}[theorem]{Definition}
   
\theoremstyle{remark}    
  \newtheorem{remark}[theorem]{Remark}
  \newtheorem{example}[theorem]{Example}

\theoremstyle{plain}
    \newtheorem{maintheorem}{Theorem}
    
    \newtheorem{maincoro}[maintheorem]{Corollary}

\usepackage{subfiles}
\usepackage{xr}

\usepackage{natbib}
\usepackage[colorlinks=true,
    citecolor=blue,
    linkcolor=black,
    urlcolor=blue
]{hyperref}

\newcommand{\nc}{\newcommand}

\nc{\RR}{\mathbb{R}}
\nc{\CC}{\mathbb{C}}
\nc{\ZZ}{\mathbb{Z}}
\nc{\QQ}{\mathbb{Q}}
\nc{\NN}{\mathbb{N}}

\nc{\Ric}{\operatorname{Ric}}
\nc{\ric}{\operatorname{ric}}
\nc{\Rm}{\operatorname{Rm}}
\nc{\adM}{\mathrm{ad}M}

\nc{\diam}{\operatorname{diam}}
\nc{\inj}{\operatorname{inj}}
\nc{\Vol}{\operatorname{Vol}}

\title{Ricci flows and bracket flows on twisted principal bundles}
\author{Henry Wood}
\date{}

\begin{document}
   \raggedbottom
\begin{abstract}
    We show that Ricci flows of invariant metrics on twisted principal bundles with nilpotent fibres over \(S^1\) are immortal and type III. 
    After normalisation and lift to the universal cover, a sequence of metrics converges in the Cheeger–Gromov topology to a homogeneous expanding Ricci soliton. 
    The proof extends the bracket flow technique from homogeneous Ricci flows to invariant Ricci flows on twisted principal bundles. 
\end{abstract}
\maketitle

\setcounter{tocdepth}{1}
\tableofcontents

\section{Introduction}

A \textit{Ricci flow} is a family of Riemannian metrics \(g(t)\) on a smooth manifold that solves
\begin{equation}\label{equ:Ricci-flow}
    \frac{\partial g}{\partial t}  = - 2\ric_g, \qquad g(0) = g_0.
\end{equation}
A central problem is to determine whether the Ricci flow encounters a finite-time singularity or admits a long-time solution. 
Solutions to \eqref{equ:Ricci-flow} that exist for all \(t>0\) are called \textit{immortal} Ricci flows.
Key to the study of immortal Ricci flows is understanding their long-time asymptotics. 
In certain cases, immortal Ricci flows converge to expanding Ricci solitons, self-similar solutions to the Ricci flow that evolve only by homothetic expansion and pullback by diffeomorphisms \cite{lott2007long, lott2010dim, bohm2015long, bohmlaf2018, lott2026ricciflow3torusgoes}. 
\newline
\indent
The long-time behaviour of immortal Ricci flows on closed \(3\)-manifolds admitting locally homogeneous metrics modelled on \(\RR^3\), \({\rm Nil}\), and \({\rm Sol}\) is now characterised. 
Up to finite covers, these manifolds are torus bundles over \(S^1\) with local \(\mathsf{T}^2\)-actions related by group automorphisms prescribed by the monodromy.
For each of the three geometries, the locally homogeneous Ricci flows on these manifolds were described in \cite{isenberg1992ricci}. 
Lott--{\v{S}}e{\v{s}}um prove that if the initial metric is invariant under the local \(\mathsf{T}^2\)-actions, then its Ricci flow is immortal and satisfies type III curvature estimates \cite{lottsesum2014}. 
Recent work by Lott shows that, for arbitrary immortal Ricci flows on these manifolds without symmetry, the blowdowns of the metrics lifted to the universal cover converge to homogeneous expanding Ricci solitons \cite{lott2026ricciflow3torusgoes}. 
\newline
\indent
These torus bundles belong to a broader class of manifolds with local symmetries, called \textit{twisted principal bundles}: fibre bundles that have Lie group fibres and transition functions consisting of smoothly varying left translations and locally constant group automorphisms.
The significance of this class of manifolds is illustrated by Lafuente--Thompson's recent constructions of expanding Ricci solitons with constant scalar curvature.
Starting in dimension four, their constructions yield many novel examples that are not locally homogeneous \cite{lafuente2025expandingriccisolitonshiggs}.
In a wider context, this type of local symmetry plays a key role in the study of Ricci flows that collapse with bounded curvature \cite{lott2010dim, cheeger1992nilpotent}.
\newline
\indent
Our first main theorem establishes long-time existence, together with type III estimates, for Ricci flows of invariant metrics on twisted principal bundles with nilpotent fibres over \(S^1\) in arbitrary dimension.

\begin{maintheorem}\label{thm:A}
    Let \(\mathsf{N}\) be a nilpotent Lie group.
    Ricci flows of invariant metrics \(g(t)\) on twisted principal \(\mathsf{N}\)-bundles over \(S^1\) are immortal and satisfy the type III estimates
    \[
    |\Rm_{g(t)}|_{g(t)} \leq \frac{C}{t}, \qquad \diam(S^1, \check{g}(t)) \leq C_0\sqrt{1+t}, 
    \]
    where \(\check{g}(t)\) is the quotient metric. 
    Here \(C >0\) is a constant depending only on dimension, and \(C_0 >0\) is a constant depending only on the initial metric. 
\end{maintheorem}

\autoref{thm:A} extends the long-time existence results and curvature estimates proved by Lauret for Ricci flows of left-invariant metrics on nilpotent Lie groups \cite{lauret2011nilmanifolds} and by Gindi--Streets for Ricci flows of invariant metrics on four-dimensional nilpotent principal bundles over \(S^1\) \cite{gindi2023four}.  
Lauret's proof relied in part on a technique called \textit{bracket flow}, an equivalent formulation of homogeneous Ricci flow in which the Lie bracket evolves while the inner product on the Lie algebra remains fixed. 
The bracket flow was first developed for nilpotent Lie groups by Guzhvina \cite{guzhvina2008ricci} and later by Lauret for simply-connected nilmanifolds and homogeneous spaces \cite{lauret2011nilmanifolds, lauret2013RicciFlowOf}.
Subsequent applications include Lafuente's scalar curvature blowup criterion for homogeneous Ricci flows \cite{lafuente2015scalar} and B\"ohm--Lafuente's proof that every sequence of blowdowns of an immortal homogeneous Ricci flow subconverges to a homogeneous expanding Ricci soliton \cite{bohmlaf2018}, with further applications to other geometric flows \cite{lauret2016geometric,thompson2024bach}.  
\newline
\indent
To prove \autoref{thm:A}, we extend the homogeneous bracket flow to twisted principal bundles.
For simplicity, let us explain our construction in the particular case of Ricci flows \((M, g)\) that are invariant under a free and proper isometric action of a connected Lie group \(\mathsf{G}\) with \emph{arbitrary cohomogeneity}.
The \(\mathsf{G}\)-invariant vector fields correspond to sections of a vector bundle \(TM/\mathsf{G}\to M/\mathsf{G}\). 
This bundle is endowed with a Lie bracket on sections satisfying some natural axioms, which in this particular case is simply the bracket of \(\mathsf{G}\)-invariant vector fields on \(M\).
Thus, \(TM/\mathsf{G}\) is a \textit{Lie algebroid}.
\newline
\indent
Sections of the \textit{adjoint bundle}, the Lie algebra bundle associated to \(M \to M/\mathsf{G}\) by the adjoint representation \(\adM = M \times_{\rm Ad} \mathfrak{g}\), correspond to invariant vertical vector fields. 
The differential of the submersion \(M \to M/\mathsf{G}\) induces the short exact sequence of vector bundles, 
\begin{equation}\label{equ:intro-exact-sequence}
    0 \to \adM \to TM/\mathsf{G} \to T(M/\mathsf{G}) \to 0. 
\end{equation}
Splittings of this sequence \eqref{equ:intro-exact-sequence} correspond to principal connections on \(M\to M/\mathsf{G}\). 
The \(\mathsf{G}\)-invariant metrics on \(M\) are in one-to-one correspondence with bundle metrics on \(TM/\mathsf{G}\), which are determined by a bundle metric on \(\adM\), a Riemannian metric on \(M/\mathsf{G}\), and a principal connection. 
The main idea for the higher-cohomogeneity bracket flow is to evolve the Lie algebroid bracket of \(TM/\mathsf{G}\), together with the Riemannian metric on \(M/\mathsf{G}\), and the connection, while keeping the bundle metric on \(\adM\) fixed.  
When \(M/\mathsf{G}\) is a point, \(TM/\mathsf{G}\) is the Lie algebra \(\mathfrak{g}\) and we recover Lauret's homogeneous bracket flow. 
\newline
\indent
Our next theorem establishes equivalence between the Ricci flow of invariant metrics and the bracket flow (see \autoref{thm:precise-B} for a more precise statement).

\begin{maintheorem}\label{thm:B}
    On twisted principal bundles, the Ricci flow of invariant metrics and the bracket flow are equivalent as deformations of Riemannian Lie algebroids. 
\end{maintheorem}

A key ingredient in the proof is Uhlenbeck's trick of moving frames.
Given an invariant Ricci flow, take a moving frame of invariant vertical vector fields that remains orthonormal along the flow. 
We then pull back the fixed Lie algebroid bracket by this change of frame, to obtain a curve of Lie algebroid brackets. 
Given either the Ricci flow or the bracket flow, the other is obtained by solving for a time-dependent change of frame determined by a fibrewise ordinary differential equation. 
Consequently, both flows have the same maximal time of existence, and curvature estimates obtained for the bracket flow hold under the Ricci flow.
\newline
\indent
To understand the convergence and long-time behaviour of the Ricci flows in \autoref{thm:A}, we consider the normalised Ricci flow:
\begin{equation}\label{equ:normalised-Ricci-flow}
    \frac{\partial g}{\partial t}  = - 2\ric_g - g, \qquad g(0) = g_0. 
\end{equation}
The two flows are related in the following sense: if \(g(t)\) is a normalised Ricci flow \eqref{equ:normalised-Ricci-flow}, then \(\bar{g}(s) = sg(\log s)\) is a solution to the unnormalised Ricci flow \eqref{equ:Ricci-flow}.
Hence, immortal Ricci flows defined for all \(s>0\) correspond to eternal normalised Ricci flows, that is, solutions to \eqref{equ:normalised-Ricci-flow} that exist for all times. 
In this correspondence, \(g(t) = e^{-t}\bar{g}(e^t)\) and \(g_0 = \bar{g}(1)\). 
The type III estimates of \autoref{thm:A} imply that the corresponding normalised Ricci flows have uniformly bounded Riemannian curvature and quotient diameter for all positive time.  
For this reason, we study the time translates \(g_k(t) := g(t +t_k)\) of the normalised Ricci flow. 
If \(s_k = e^{t_k}\), then \(s_k^{-1}\bar{g}(s_k s) = s g_k(\log s)\), so the time translates describe the blowdowns of the unnormalised Ricci flow. 
Indeed, if \(g_k(t)\) converges to \(g_\infty(t)\), then the blowdowns converge to \(sg_\infty(\log s)\). 
\newline
\indent
Gindi derives rigidity results for blowdown limits of the invariant Ricci flow on nilpotent principal bundles, assuming integrability of the horizontal distribution along the flow \cite[Thm. 8.1]{gindi2026curvature}. 
For Heisenberg fibres, the limits are locally expanding Ricci solitons \cite[Thm. 10.3]{gindi2026curvature}. 
For one-dimensional base, the horizontal distribution always remains integrable along the flow. 
\newline
\indent
Our final theorem demonstrates that the normalised Ricci flow of invariant metrics on twisted principal bundles with nilpotent fibres over \(S^1\), once lifted to the universal cover, subconverges to a solvsoliton.

\begin{maintheorem}\label{thm:C}
    Let \(\mathsf{N}\) be a simply-connected nilpotent Lie group.
    Let \(M \to S^1\) be a twisted principal \(\mathsf{N}\)-bundle. 
    Let \(g(t)\) be a normalised Ricci flow of invariant metrics on \(M\).  
    Let \(\tilde{g}(t)\) be the lift of \(g(t)\) to the universal cover \(\tilde{M}\). 
    Fix a point \(\tilde{p} \in \tilde{M}\).
    Then, there exists a sequence \(t_k \to \infty\) and a simply-connected solvable Lie group \(\mathsf{S}\) equipped with a solvsoliton metric \(g_\infty\) such that 
    \[
    (\tilde{M}, \tilde{g}(t_{k}), \tilde{p}) \underset{k \to \infty}{\longrightarrow} (\mathsf{S}, g_\infty, \mathsf{e})
    \]
    in the smooth pointed Cheeger--Gromov topology. 
\end{maintheorem}

Recall that a solvsoliton is an invariant expanding Ricci soliton on a simply-connected solvable Lie group that evolves under Ricci flow by homothetic expansion and pullback by Lie group automorphisms \cite{lauret2011solvmanifolds}. 
To prove \autoref{thm:C}, we first prove a compactness theorem for normalised Ricci flows of invariant metrics on twisted principal bundles over \(S^1\) with simply-connected nilpotent fibres.  
The compactness argument combines the type III estimates of \autoref{thm:A} with Cheeger--Gromov--Taylor's injectivity radius estimate \cite[Thm. 4.7]{cheegergromovtaylor1982finite}, Hamilton's compactness theorem \cite{hamilton1996compactness}, and Heber's compactness of Killing fields \cite[p. 330]{heber1998noncompact}. 
The assumptions that the fibres are nilpotent and simply-connected provide the volume non-collapsing needed for the Cheeger--Gromov--Taylor injectivity radius estimate.
If the diameter of the base collapses, then the limit of the universal covers is homogeneous. 
We appeal to \cite{bohmlaf2018, bohmlafuente2023noncompact} to obtain a solvsoliton limit.  
If the base does not collapse, then  we can derive properties of the limit using Gindi--Streets monotone energy functional \cite[Def. 4.1]{gindi2023four}:
\begin{equation}\label{equ:Gind-Streets-energy-intro}
    \mathcal{I}(g) = \int_{S^1} |L|^2 + \frac{1}{2}\, ds. 
\end{equation}
After obtaining a limit Ricci flow, we take the corresponding bracket flow.
The properties obtained from \eqref{equ:Gind-Streets-energy-intro} reduce this bracket flow to a homogeneous bracket flow on the fibres, with a fixed connection and quotient metric. 
We integrate a monotone quantity for homogeneous normalised bracket flow \cite{lauret2011nilmanifolds}:
\begin{equation}
     \mathcal{F}(\mu, h, \check{g}) = \int_{S^1} \hat{R}_\mu + |\hat{\Ric}_\mu|^2\, ds.
\end{equation}
The monotonicity of \(\mathcal{F}\) shows that a further subsequential limit has nilsoliton fibres, while retaining the properties established by \eqref{equ:Gind-Streets-energy-intro}.  
Together, these conditions imply that the limit metric is a solvsoliton, proving \autoref{thm:C}. 
\newline
\indent
For twisted principal bundles with unipotent monodromy, we use the limit properties from the proof of \autoref{thm:C} to show collapse of the quotient \(S^1\). 
Arguing by contradiction, we assume that the quotient diameter remains bounded below along a sequence of times of the normalised Ricci flow. 
The corresponding limit has unipotent holonomy and a parallel non-zero shape operator. 
These properties are contradictory, proving the following corollary (see \autoref{coro:D-precise} for a more precise statement).

\begin{maincoro}\label{coro:D}
    Let \(\mathsf{N}\) be a simply-connected nilpotent Lie group acting with cohomogeneity one on \(M\), with orbit space \(M/\mathsf{N} \cong S^1\).  
    If \(g(t)\) is a normalised Ricci flow of \(\mathsf{N}\)-invariant metrics on \(M\), then
    \[
    \lim_{t \to \infty}\diam(M/\mathsf{N} , \check{g}(t)) =  0,
    \]
    where \(\check{g}(t)\) is the quotient metric.
\end{maincoro}

For principal \(\mathsf{N}\)-bundles over \(S^1\), the collapse of the base given in \autoref{coro:D} is incompatible with the convergence through principal bundle isomorphisms assumed in \cite[Thm. 1.16, Thm. 11.1]{gindi2026curvature}. 
As observed in \cite[\S 1.1.5]{gindi2026curvature}, these results can be extended to twisted principal bundles, where monodromy can prevent the base from collapsing. 

\begin{example}
    Lott--{\v{S}}e{\v{s}}um provide the following example \cite[eq. 3.9]{lottsesum2014}. 
    Take \(M\) to be the twisted principal \(\RR^2\)-bundle over \(S^1\) with hyperbolic monodromy given by 
    \[
    H = \begin{pmatrix}
        2 & 1 \\ 1 & 1
    \end{pmatrix}.
    \]
    Let \(x_+, x_- \in \RR^2\) be eigenvectors of \(H\) with corresponding eigenvalues \(\lambda = \frac{1}{2}(3 +\sqrt{5})\) and \(\lambda^{-1} = \frac{1}{2}(3 -\sqrt{5})\). 
    Let \(c = \log \lambda\). 
    Then 
    \[
    g(t) = 4c^2 dy^2 + e^{-t} (e^{2cy} dx_+^2  + e^{-2cy} dx_-^2)
    \]
    is an invariant normalised Ricci flow on \(M \to S^1\) with a constant quotient metric \mbox{\(\check{g}(t) =  4c^2 dy^2\)}. 
\end{example}

\subsection{Acknowledgements}
The author would like to thank his supervisors, Ramiro Lafuente and Kyle Broder, for their guidance, helpful discussions, and feedback on earlier versions of this work.  
The author gratefully acknowledges access to OpenAI’s ChatGPT for Academic Researchers program through a research workspace awarded to Kyle Broder.
This work was supported by an Australian Government Research Training Program (RTP) Scholarship.

\subsection{AI disclosure}
No AI tools were used in the proofs of \autoref{thm:A} or \autoref{thm:B}. ChatGPT 5.6 Sol, 6 Astra, and 6.1 Sol assisted in the proofs of \autoref{thm:C} and \autoref{coro:D}. 
The author takes full responsibility for the mathematical content of this paper. 
No text was written by AI.

\section{Preliminaries}

\subsection{Lauret's bracket flow and Uhlenbeck's trick}\label{subsec:Lauret's}
To begin, we recall Lauret's bracket flow on Lie groups
and demonstrate its equivalence to the Ricci flow ~\cite{lauret2011nilmanifolds, lauret2013RicciFlowOf}. 
Let \(\mathsf{G}\) be a connected Lie group and let \(g(t): t\in [0, T)\) be a Ricci flow of left-invariant metrics with initial metric \(g(0)=g_0\).
\newline
\indent
By a trick attributed to Karen Uhlenbeck, we simplify the Ricci flow equation by applying an appropriate change of frame (cf. \cite[p.~155]{hamilton1986four}).
At \(t=0\), choose a left-invariant \(g_0\)-orthonormal frame \(\{X_i\}\). 
Evolve this frame by 
\[
\partial_t X_i = \Ric_{g(t)} (X_i),
\]
where \(\Ric_{g(t)}\) is the Ricci endomorphism. 
This provides a moving frame \(\{X_i(t)\}\) on \(\mathsf{G}\) which remains \(g(t)\)-orthonormal under Ricci flow:
\begin{align*}
    \pdv{}{t} g(X_i, X_j) &=
    \left(\partial_t g  \right)(X_i, X_j) +
    g\left(\partial_t X_i, X_j\right)
    + g\left(X_i, \partial_tX_j\right)
    \\
    &= -2\ric_{g(t)}(X_i, X_j) + g(\Ric_{g(t)}(X_i), X_j) + g(X_i, \Ric_{g(t)}(X_j)) = 0.
\end{align*}
The structure constants are \(\mu_{ij}^k(t)X_k(t) = [X_i(t), X_j(t)]\). 
By differentiating this equation, we get the ODE for the structure constants
\[
\frac{d\mu_{ij}^k}{dt}
X_k 
= [\Ric_{g(t)} (X_i), X_j] + [ X_i, \Ric_{g(t)} (X_j)]
-\Ric_{g(t)}[X_i, X_j].
\]
To express this evolution invariantly, let  \(\langle \cdot, \cdot\rangle_0\) be the inner product on \(\mathfrak{g}\) given by evaluating \(g_0\) at the identity. 
Take \(\mu \in \Lambda^2\mathfrak{g}^* \otimes \mathfrak{g}\) such that \((\mathfrak{g}, \mu)\) and \((\mathfrak{g}, [\cdot, \cdot])\) are isomorphic as Lie algebras. 
By the third theorem of Lie, there exists a connected and simply-connected Lie group \(\tilde{\mathsf{G}}_\mu\) with Lie algebra \((\mathfrak{g}, \mu)\).
Let \(\tilde{g}_\mu\) be the left-invariant metric on \(\tilde{\mathsf{G}}_\mu\) such that \(\tilde{g}_\mu (\mathsf{e}) = \langle \cdot, \cdot \rangle_0\).  
Let \(\Ric_\mu\) be the Ricci endomorphism of \((\tilde{\mathsf{G}}_\mu, \tilde{g}_\mu)\) at the identity.
Consider
\begin{equation}\label{equ:bracket-flow}
    \frac{d\mu}{dt} = \delta_\mu(\Ric_\mu),
    \qquad \mu(0) = [\cdot, \cdot]
\end{equation}
where \(\delta_\mu :{\rm End}(\mathfrak{g}) \to \Lambda^2 \mathfrak{g}^* \otimes \mathfrak{g}\) is given by
\[
\delta_\mu(\alpha) :=\frac{d}{dt}\Big|_{t = 0} e^{-\alpha t}\mu(e^{\alpha t} \cdot, e^{\alpha t}\cdot) =\mu(\alpha \cdot, \cdot) + \mu(\cdot, \alpha \cdot) - \alpha \mu(\cdot, \cdot),
\]
for all \(\alpha \in {\rm End}(\mathfrak{g})\). 
The ordinary differential equation \eqref{equ:bracket-flow} is called the \textit{bracket flow}.
As the time derivative \(\tfrac{d\mu}{dt}\) is tangent to the orbit \(\mathsf{GL}(\mathfrak{g})\cdot [\cdot, \cdot]\), there exists a unique smooth one-parameter family of linear transformations \(\Phi(t)\) satisfying \(\mu(t) = \Phi(t) \cdot [\cdot, \cdot]\) and solving
\[
\frac{d\Phi}{dt} = - \Ric_\mu \Phi, \qquad \Phi(0) = \operatorname{Id}.  
\]
The maps \(\Phi(t) : (\mathfrak{g}, [\cdot, \cdot]) \to (\mathfrak{g}, \mu(t))\) are Lie algebra isomorphisms. 
They integrate uniquely to Lie group isomorphisms \(\tilde{\varphi}(t): \tilde{\mathsf{G}} \to \tilde{\mathsf{G}}_{\mu(t)}\) with \(d\tilde{\varphi}(t)_\mathsf{e} = \Phi(t)\). 
Now, if \(\Gamma\) is the kernel of the universal covering homomorphism, then define \(\mathsf{G}_{\mu(t)} = \tilde{\mathsf{G}}_{\mu(t)}/\tilde{\varphi}(t)(\Gamma)\).
The metrics \(\tilde{g}_{\mu(t)}\) descend to \(g_{\mu(t)}\), and the Lie group isomorphisms \(\tilde\varphi(t)\) descend to \({\varphi}(t): \mathsf{G} \to \mathsf{G}_{\mu(t)}\). 
By pulling back \(g_{\mu(t)}\) by the Lie group isomorphisms \(\varphi(t): \mathsf{G} \to \mathsf{G}_{\mu(t)}\), we obtain the Ricci flow of left-invariant metrics \(g(t)\):
\begin{align*}
    \frac{d}{dt} (\varphi(t)^*g_{\mu(t)})(\mathsf{e}) &= \frac{d}{dt} \langle \Phi(t)\cdot, \Phi(t)\cdot\rangle_0 
    \\
    &= \langle\Phi'(t) \cdot, \Phi(t)\cdot\rangle_0 + \langle\Phi(t) \cdot, \Phi'(t)\cdot\rangle_0
    \\
    &= -\langle \Ric_{\mu(t)} \Phi(t) \cdot, \Phi(t)\cdot\rangle_0 - \langle\Phi(t) \cdot,  \Ric_{\mu(t)} \Phi(t)\cdot\rangle_0
    \\
    &=-2\langle\Ric_{\mu(t)}\Phi(t) \cdot, \Phi(t)\cdot\rangle_0
    \\
    &= -2\ric_{\varphi(t)^*g_{\mu(t)}}(\mathsf{e)}.
\end{align*}
The Ricci flow \(g(t)\) and the bracket flow \(g_{\mu(t)}\) differ only by pullback by time-dependent diffeomorphisms \cite[Thm. 3.3]{lauret2013RicciFlowOf}. 
The curvatures and all Riemannian invariants of the Ricci flow can be studied along the bracket flow.

\subsection{Riemannian geometry on twisted principal bundles}\label{subsec:Riemann-Geometry}
Here, we define twisted principal bundles, the class of manifolds considered in our main theorems.  
The material in this subsection is based on   
\cite[\S 2]{lafuente2025expandingriccisolitonshiggs} and
\cite[\S 3.1]{lottsesum2014}.

\begin{definition}
    Let \(\mathsf{G}\) be a Lie group. 
    A \textnormal{twisted principal} \(\mathsf{G}\)-\textnormal{bundle} is a fibre bundle \(\pi: M \to B\) with fibre \(\mathsf{G}\) and structure group contained in \({\rm Aut}(\mathsf{G}) \ltimes \mathsf{G}_L\) where the topology on \({\rm Aut}(\mathsf{G})\) is discrete. 
    Here, \(\mathsf{G}_L\) is the group of left-multiplications of \(\mathsf{G}\) acting  on itself. 
\end{definition}

The structure group of a twisted principal bundle need not commute with right group multiplication \(\mathsf{G}_R\), and so the fibrewise \(\mathsf{G}\)-actions do not, in general, extend to a global action on \(M\).  
Let \((U_\alpha, \psi_\alpha)\)
be a local trivialisation of \(M\), and define the right \(\mathsf{G}\)-action on \(M_{U_\alpha} := \pi^{-1}(U_\alpha)\) by
\begin{equation}\label{equ:local-G-actions}
R^{\alpha} : M_{U_\alpha} \times \mathsf{G} \to M_{U_\alpha}, \qquad  R^{\alpha}_\mathsf{g}(p) = \psi_\alpha^{-1}( \psi_\alpha(p) \cdot\mathsf{g}). 
\end{equation}
If \((U_\beta, \psi_\beta)\) is an overlapping local trivialisation with transition function \(\tau_{\alpha\beta}(x) = \varphi_{\alpha\beta} \circ L_{\mathsf{h}(x)}\) where \(\varphi_{\alpha\beta}\) is constant, 
then 
\[
R^\alpha_\mathsf{g} (p)
= \psi_\beta^{-1}( \psi_\beta(p) \cdot\varphi_{\alpha \beta}(\mathsf{g})) = R^\beta_{\varphi_{\alpha\beta}(\mathsf{g})} (p).
\]

\begin{example}[M\"obius band]
    The M\"obius band is a twisted principal \(\RR\)-bundle over \(S^1\) with structure group \(\{\pm \operatorname{Id}\} \ltimes \RR\). 
    Cover \(S^1\) by two trivialisations \((U_\alpha, \psi_\alpha)\) and \((U_\beta, \psi_\beta)\) with transition functions \(\operatorname{Id}\) and \(-\operatorname{Id}\). 
    Then, on one of the overlaps (connected components of \(U_\alpha \cap U_\beta\)), 
    \[
    R^\alpha_{\mathsf{g}}(p) = \psi_\alpha^{-1}( \psi_\alpha(p) + \mathsf{g}) 
    = \psi_\beta^{-1}( \psi_\beta(p) -\mathsf{g}) 
    = R^\beta_{-\mathsf{g}} (p).
    \]
\end{example}

\begin{example}[Torus bundles]
    Let \(H \in \mathsf{SL}_2(\ZZ)\). 
    The following is a twisted principal \(\mathsf{T}^2\)-bundle over \(S^1\), 
    \[
    M = (\RR \times \mathsf{T}^2)/\sim, \qquad (s +1, x) \sim (s, Hx).
    \]
    If \(H\) is elliptic, i.e., \(H\) is of finite order, then \(M\) is finitely covered by \(T^3\), so \(M\) admits a flat geometry. 
    If \(H\) is parabolic, that is, \(|\tr H| = 2\), but \(H\) is not elliptic, then \(M\) admits a \({\rm Nil}\) geometry. 
    If \(H\) is hyperbolic, that is, \(|\tr H| >2\), then \(M\) admits a \({\rm Sol}\) geometry. 
\end{example}

Although there may not be a globally defined \(\mathsf{G}\)-action on \(M\), the local actions \eqref{equ:local-G-actions} provide well-defined natural transformation groups \(\mathsf{G}_U \subset {\rm Diff}(M_{U})\) that act freely and transitively on the fibres. 
A choice of local trivialisation gives an isomorphism \(\mathsf{G} \simeq \mathsf{G}_{U}\). 
In lieu of a global action, these natural transformation groups are sufficient to provide a well-defined notion of \(\mathsf{G}\)-invariance on \(M\). 
\newline
\indent
A tensor on a twisted principal bundle \(\pi: M \to B\) is called \textit{invariant} if it is invariant under the actions of the local transformation groups \(\mathsf{G}_{U}\) for all trivialising neighbourhoods \(M_U \simeq U\times \mathsf{G}\). 
The bundle of invariant vertical vector fields is a Lie algebra bundle called the \textit{adjoint bundle} \(\adM \to B\). 
If \(\{U_\alpha \times \mathsf{G}\}\) is a bundle atlas of \(M\), then \(\{U_\alpha \times \mathfrak{g}\}\) is a bundle atlas of \(\adM\) with transition functions given by the group homomorphism 
\[
\operatorname{Aut}(\mathsf{G})\ltimes \mathsf{G}_L 
\to \operatorname{Aut}(\mathfrak{g}), 
\qquad 
(\varphi, \mathsf{g}) \mapsto d\varphi_e  \circ \operatorname{Ad}(\mathsf{g}).
\]
As the structure group of \(\adM\) is contained in \({\rm Aut}(\mathfrak{g})\), the Lie bracket on \(\mathfrak{g}\) induces a canonical section \([\cdot, \cdot] \in \Gamma(\Lambda^2\adM^* \otimes \adM)\) whose restrictions to fibres are Lie brackets. 
\newline
\indent
A \textit{horizontal distribution} \(\mathcal{H}\subset{ TM}\) on a twisted principal bundle is a distribution that is complementary to the vertical distribution \(\mathcal{V} = \ker \pi_*\). 
An invariant horizontal distribution defines a \textit{principal connection}.
A principal connection is a vector bundle map
\[
\omega:TM \to \pi^*\adM   
\]
which when restricted to the vertical distribution is the canonical identification of \(\mathcal{V}_p\) with \(\adM_{\pi(p)}\). 
Vice versa, the kernel of a principal connection is a given invariant horizontal distribution. 
A horizontal lift of \(X \in \mathfrak{X}(B)\) is denoted \(\bar{X}\).
Given an invariant vertical vector field \(U\), denote the corresponding section of the adjoint bundle as \(U_{\rm ad} := \omega(U)\).
\newline
\indent
Equip \(M \to B\) with an invariant Riemannian metric \(g\). 
The orthogonal complement of the vertical distribution determines a principal connection \(\omega\), the restriction to invariant vertical vector fields determines a vector bundle metric \(h\) on \(\adM\), and the metric on the horizontal distribution descends to a Riemannian metric \(\check{g}\) on \(B\). 
The invariant Riemannian metric uniquely decomposes as 
\[
g = \omega^*h + \pi^*\check{g} =h(\omega\cdot, \omega\cdot) +\check{g}(\pi_*\cdot, \pi_*\cdot). 
\]
Conversely, every such triple \((h, \omega, \check{g})\) uniquely determines an invariant Riemannian metric on \(M\). 
There is a bijection
\begin{equation}\label{equ:reduce-data-bij}
    \{\text{invariant Riemannian metrics \(g\) on \(M\)}\} \longleftrightarrow
    \left\{
    \begin{array}{cc}
     \text{vector bundle metrics \(h\) on } \adM  
     \\
     \text{principal connections \(\omega\)}
     \\
    \text{Riemannian metrics  \(\check{g}\) on } B
    \end{array}
    \right\}.
\end{equation}
\indent
Let \(U, V, W\) be vertical vector fields on \(M\), and let \(X, Y, Z\) be vector fields on \(B\). 
If \(U\) is invariant, denote \(U_{\rm ad} \in \Gamma(\adM)\) as the corresponding section. 
A principal connection on \(M \to B\) determines horizontal lifts \(\bar{X}\) that are invariant under the local natural transformations \(\mathsf{G}_U\). 
This induces an affine connection on \(\adM\):
\begin{equation}\label{equ:def-principal-affine-connection}
    \nabla^{\rm ad}: \mathfrak{X}(B) \times \Gamma(\adM) \to \Gamma(\adM), 
\qquad
\nabla^{\rm ad}_X U_{\rm ad} := [\bar{X}, U]_{\rm ad}.
\end{equation}
\indent
Equip \(M\) with an invariant Riemannian metric \(g\). 
Let \(\nabla^M\) be the Levi-Civita connection of \((M, g)\), and define 
\[
\nabla: \mathfrak{X}(B) \times \Gamma(\adM) \to \Gamma(\adM), 
\qquad \nabla_X U_{\rm ad} := (\mathcal{V}\nabla^M_{\bar{X}} U)_{\rm ad},
\]
where \(\mathcal{V}\) denotes the vertical component of the vector field.
This connection is metric on \((\adM, h)\) and is characterised by Koszul's formula:
\[
h(\nabla_{X} U_{\rm ad}, V_{\rm ad}) =\frac{1}{2} \left(Xh(U_{\rm ad}, V_{\rm ad}) + h(\nabla^{\rm ad}_X U_{\rm ad}, V_{\rm ad})- h(U_{\rm ad}, \nabla^{\rm ad}_X V_{\rm ad} )\right). 
\]
As such, the metric connection \(\nabla\) may be defined intrinsically on \((\adM, h)\) without the need to compute the Levi-Civita connection of \((M, g)\). 
\newline
\indent
The space of connections on a vector bundle is affine and modelled on the endomorphism-valued one-forms. 
There exists an \(L \in \Omega^1(B;{\rm End}(\adM))\) such that 
\[
\nabla^{\rm ad} = \nabla - L. 
\]
By the torsion freeness of \(\nabla^M\), this one-form is the shape operator of the fibres: 
\[
L_{{X}} U_{\rm ad} = (\mathcal{V} \nabla^{M}_{U} \bar{X})_{\rm ad}. 
\]
If \(\{U_i\}\) is a local orthonormal vertical frame on \((M, g)\) and \(\{X_j\}\) is a local orthonormal frame on \((B, \check{g})\), 
then the \textit{mean curvature vector} of the fibres is given by 
\[
N = \sum_i \mathcal{H}\nabla^M_{U_i} U_i = -\sum_j (\tr L_{X_j})\bar{X}_j. 
\]
This is a basic horizontal vector field in \(TM\). 
We also write \(N\) for its projection \(\pi_*N\). 
The \textit{homogeneous mean curvature vector} of the fibres is characterised by 
\[
h(H, U_{\rm ad}) = {\rm tr}(  \operatorname{ad}_{U_{\rm ad}}). 
\]
Let \(\hat\nabla\) be the Levi-Civita connection of the fibres of \((M, g) \to (B, \check{g})\) given by \(\hat{\nabla}_U V = \mathcal{V}\nabla^M_U V\). 
Its restriction to invariant vertical vector fields is again invariant, thus providing the bundle map 
\[
\hat{\nabla}: \adM \otimes \adM  \to \adM, 
\qquad \hat\nabla_{U_{\rm ad}} V_{\rm ad} := (\mathcal{V} \nabla_{U}^M V)_{\rm ad}. 
\]
Let \(\check{\nabla}\) be the Levi-Civita connection of the base \((B, \check{g})\). 
The curvature of the principal connection \(\omega\) is given by 
\[
F \in \Omega^2(B; \adM), \qquad
F_{X, Y} := -\omega[ \bar X, \bar Y]. 
\]
O'Neill's integrability tensor, viewed as an \(\adM\)-valued two-form on \(B\), is given by 
\[
A \in \Omega^2(B; \adM), \qquad 
A_{X, Y} : = \frac{1}{2} \omega[\bar{X}, \bar{Y}]= -\frac{1}{2}F_{X, Y}.
\]
For \(X, Y \in \mathfrak{X}(B)\) and \(U_{\rm ad} \in \Gamma(\adM)\), write
\[
A_X(Y) := A_{X, Y}, \qquad (A^*U_{\rm ad})(X, Y) := h(A_{X, Y}, U_{\rm ad}), \qquad (LU_{\rm ad})(X) := L_X U_{\rm ad}.
\]
We use \(g\) for the induced tensor inner products, obtained by full contraction with \(h\) and \(\check{g}\), and use the corresponding tensor norms. 
\newline
\indent
Next, we recall O’Neill’s Ricci curvature formulas, expressed in the notation and tensor inner product conventions of this subsection \cite[9.36]{besse2007einstein}. 
In the following, we drop the subscript \({\rm ad}\) for sections of the adjoint bundle. 

\begin{prop}[O’Neill’s Ricci curvature formulas]
    Let \(M \to B\) be a twisted principal bundle equipped with an invariant metric \(g\). 
    Then, 
    \begin{align}
        \ric_g(U, U) 
        &= \hat \ric_g(U, U) + |A^*U|^2 + h(L_N U, U)+ h((\nabla^* L) U, U),
        \\
        \ric_g(U, \bar X)&=
        -h({\hat\nabla}^* L_X, U) -  h((\nabla^{*} A)X, U) - g(A_X, L U) + h(A_{X, N}, U),
        \\
        \ric_g(\bar X, \bar X)&= \ric_{\check{g}}(X, X) - 2|A_X|^2 -|L_X|^2 + \check{g}(\check\nabla_X N,X).
    \end{align}
    Here \(\hat\ric_g\) is the Ricci tensor of the fibres, \(A^*
    \) is the adjoint of \(A\) with respect to \(h\), and 
    \(\nabla^*, \hat{\nabla}^*\) are the formal adjoints to the connections \(\nabla\), \(\hat{\nabla}\):
    \[
    \nabla^* \eta = -\sum_i (\nabla_{X_i} \eta)_{X_i}, \qquad
    \hat{\nabla}^* \xi =   -\sum_j (\hat{\nabla}_{U_j} \xi)_{U_j},
    \]
    for \(n\)-forms \(\eta\) on \(B\) with values in tensor powers of \(\adM\), and tensor fields \(\xi\) on \(\adM\), where \(\{X_i\}\) is a 
    local orthonormal frame on \((B, \check{g})\), and \(\{U_i\}\) is a 
    local orthonormal frame  on \((\adM, h)\). 
\end{prop}

These equations simplify considerably when the base is one-dimensional.
Notably, all one-dimensional distributions are integrable, so \(A=0\). 
To simplify notation, for a local unit vector field \(X\) on \(B\), write \(L = L_X\) and \(\nabla^*L = -(\nabla_X L)_X  = -\nabla L\). 

\begin{prop}[O’Neill’s Riemannian curvature formulas when \(\dim B = 1\)]\label{prop:Riemann-curv-dim1}
    Let \(M \to B^1\) be a twisted principal bundle equipped with an invariant metric \(g\). 
    Then, 
    \begin{align}
        \Rm_g(U, V, W, U') &= \hat\Rm_g (U, V, W, U') +h(LU, W)h(LV, U') - h(LV, W)h(LU, U'),
        \\
        \Rm_g(U, V, W, \bar{X})&= h((\hat{\nabla}_V L)U, W) -h((\hat\nabla_U L) V, W),
        \\
        \Rm_g(U, \bar{X}, \bar{X}, V) &=- h((\nabla L) U, V) - h(L^2U, V),
    \end{align}
    where \(\bar{X}\) is unit and \(\hat\Rm_g\) is the Riemannian curvature tensor of the fibres. 
\end{prop}

\begin{prop}[O’Neill’s Ricci curvature formulas when \(\dim B =1\)]\label{prop:Ricci-curv-dim1}
    Let \(M \to B^1\) be a twisted principal bundle equipped with an invariant metric \(g\). 
    Then, 
    \begin{align}
       \ric_g(U, U)
        &= \hat\ric_g(U, U) - (\tr L)h(LU, U) - h((\nabla L)U, U)
        \\
        \ric_g(\bar X, U) &=
        -h({\hat\nabla}^* L, U)
        \\
        \ric_g(\bar X, \bar X)&= -|L|^2 - \tr (\nabla L)
    \end{align}
    where \(\bar{X}\) is unit and 
    \(\hat\ric_g\) is the Ricci tensor of the fibres. 
\end{prop}

\subsection{Transitive Lie algebroids}
In this subsection, we describe the geometry of transitive Lie algebroids and their relationship to the geometry of twisted principal bundles.
For general theory, we refer to \cite{Mackenzie2005}, and for Riemannian Lie algebroids \cite{boucetta2011riemannian}.  
\newline
\indent
The vector bundle of all invariant vector fields is called the \textit{Atiyah bundle} \(TM/\mathsf{G} \to B\). 
This is perhaps misleading notation as there is no global \(\mathsf{G}\)-action on \(M\). 
Even so, this vector bundle is well-defined since the local transformation groups  \(\mathsf{G}_U \subset {\rm Diff}(M_U)\) are well-defined. 
If \(\{U_\alpha \times \mathsf{G}\}\) is a bundle atlas of \(M\), then the invariant vector fields on these trivialisations are identified by
\[
T(U_\alpha \times \mathsf{G})/\mathsf{G}
\simeq TU_\alpha \oplus (U_\alpha \times \mathfrak{g}). 
\]
The differentials of the transition functions of \(M\) descend to the local quotients and satisfy the cocycle condition, providing the transition functions of \(TM/\mathsf{G} \to B\). 

\begin{definition}
    A \textnormal{Lie algebroid} \((E, B, [\cdot, \cdot], \#)\) is a vector bundle \(E \to B\) equipped with a Lie bracket \([\cdot, \cdot]\) on the sections, with an anchor \(\#:E \to TB\). 
    An \textnormal{anchor} is a bundle map \(\#: E \to TB\) that preserves the Lie brackets, and the bracket satisfies the Leibniz identity
\[
[a, fb] = f[a, b] + \#a(f) b, 
\qquad a, b \in \Gamma(E), \quad f \in C^\infty(B).
\]
A Lie algebroid is called \textnormal{transitive} if its anchor is surjective. 
\end{definition}


\begin{example}
    As the bracket of invariant vector fields on a twisted principal \(\mathsf{G}\)-bundle \(M\to B\) is again invariant, the sections \(\Gamma(TM/\mathsf{G})\) inherit the Lie bracket of vector fields on \(M\). 
    The vector bundle \(TM/\mathsf{G}\) equipped with this Lie bracket and the anchor \(\pi_*: TM/\mathsf{G} \to TB\) is a Lie algebroid called the \textit{Atiyah algebroid} \((TM/\mathsf{G}, B, [\cdot, \cdot], \pi_*)\).
\end{example}

\begin{definition}
    Let \((E_1, B, [\cdot, \cdot]_1, \#_1)\) and \((E_2, B, [\cdot, \cdot]_2, \#_2)\) be two Lie algebroids over the same base. 
    A \textnormal{morphism} of Lie algebroids is a vector bundle map \(\Phi: E_1 \to E_2\) over the identity, such that 
    \[
    \#_2 \circ \Phi = \#_1, \qquad 
    \Phi[a, b]_1 = [\Phi(a), \Phi(b)]_2,  
    \]
    for all \(a, b \in \Gamma(E_1)\).
    An \textnormal{isomorphism} of Lie algebroids is a bijective morphism of Lie algebroids. 
\end{definition}

The kernel of the anchor \(\pi_*\) of the Atiyah algebroid consists of the invariant vertical vector fields. 
Hence, Atiyah bundles fit into the short exact sequence of vector bundles 
\begin{equation}\label{equ:Atiyah-sequence}
    0 \to \adM
\xrightarrow{\iota}
TM/\mathsf{G}
\xrightarrow{\pi_*}
TB
\to 0.
\end{equation}
A \textit{principal connection} is a splitting \(\sigma: TB \to TM/\mathsf{G}\) of the sequence \eqref{equ:Atiyah-sequence}
with \(\pi_* \circ \sigma = \operatorname{Id}_{TB}\),
or equivalently a vector bundle map \(\omega: TM/\mathsf{G} \to \adM\) with \(\omega \circ \iota = \operatorname{Id}_{\adM}\).
The two descriptions are related by
\[
\iota \circ \omega + \sigma \circ \pi_* = \operatorname{Id}_{TM/\mathsf{G}}, \qquad \omega \circ \sigma = 0.
\]
This is the decomposition of a vector field into its vertical and horizontal components:
\[
U = U^\mathcal{V} + U^\mathcal{H} = \mathcal{V}U + \mathcal{H}U  =  (\iota \circ \omega )(U) + (\sigma \circ \pi_*)(U). 
\]
The affine connection \(\nabla^{\rm ad}\) defined in \eqref{equ:def-principal-affine-connection} is the mixed term of the Lie algebroid bracket on \(TM/\mathsf{G}\):
\[
[U_{\rm ad}+ \bar{X},V_{\rm ad}+ \bar{Y} ]
=[U_{\rm ad}, V_{\rm ad}] + \nabla^{\rm ad}_X V_{\rm ad} -\nabla^{\rm ad}_Y U_{\rm ad} +[\bar{X}, \bar Y].
\]
With respect to the splitting, the last term is 
\[
[\bar{X}, \bar Y] = \sigma[X, Y] - \iota F_{X, Y} =  \overline{[X, Y]} + 2\iota A_{X, Y}. 
\]

\begin{definition}
    A \textnormal{Riemannian Lie algebroid} \((E, B, [\cdot, \cdot], \#, g)\) is a Lie algebroid with a vector bundle metric \(g\) on \(E \to B\). 
\end{definition}

\begin{definition}
    An \textnormal{isometric isomorphism of Riemannian Lie algebroids} is an isomorphism of Lie algebroids that is an isometry of metric vector bundles. 
\end{definition}

As invariant Riemannian metrics on twisted principal bundles are determined entirely by their evaluations along invariant vector fields, we have the bijection 
\begin{equation}\label{equ:Atiyah-metric-bij}
    \{\text{invariant Riemannian metrics \(g\) on \(M\)}\} \longleftrightarrow
    \{
     \text{vector bundle metrics on } TM/\mathsf{G}  
     \}.
\end{equation}

\indent
Next, we describe the local geometry of Riemannian Lie algebroids. 
A transitive Lie algebroid \((E, B, [\cdot, \cdot], \#)\) has surjective anchor, and thus there is a Lie algebra bundle \(\ker \# \to B\) and a short exact sequence
\[
    0 \to \ker \#
\xrightarrow{\iota}
E
\xrightarrow{\#}
TB
\to 0.
\]
A splitting of this sequence \(\sigma\) and its vertical projection \(\omega\) induce an affine connection on \(\ker \#\) and a curvature
\[
\nabla^{\ker }_X U_{\rm ker} = \omega[\sigma X, \iota U_{\rm ker}], 
\qquad
F_{X, Y} = -\omega[\sigma X, \sigma Y],
\qquad X \in \mathfrak{X}(B), \quad U_{\rm ker} \in \Gamma(\ker \#). 
\]
Define the deformation coboundary operator \cite[Equ. 1]{crainic2008deform}
\[
\delta_{[\cdot, \cdot]}(\alpha)(U, V)= [\alpha U, V] +  [U, \alpha V] -\alpha[U, V], \qquad \alpha \in \Gamma({\rm End}(E)),\quad  U, V \in \Gamma(E). 
\]
This operator is skew-symmetric, \(\RR\)-linear, and 
\[
\delta_{[\cdot, \cdot]}(\alpha) (U, fV) = f\delta_{[\cdot, \cdot]}(\alpha)(U, V) + (\#\alpha U)(f) V, \qquad f \in C^\infty (B). 
\]
If \(\# \circ \alpha = 0\), then \(\delta_{[\cdot, \cdot]}(\alpha)\) is \(C^\infty(B)\)-linear. 
Vector bundle metrics \(g\) on \(E\) uniquely decompose as the same metric data \((h, \omega, \check{g})\) as in subsection \ref{subsec:Riemann-Geometry}, with \(\ker\#\) instead of \(\adM\):
\[
g =  h(\omega\cdot, \omega\cdot) + \check{g}(\#\cdot, \#\cdot).
\]
Given a vector bundle metric on a transitive Lie algebroid \(E \to B\), there is a unique connection, \(\nabla: \Gamma(E) \times \Gamma(E) \to \Gamma(E)\) satisfying
\[
\nabla_{fa} b = f\nabla_a b, \qquad \nabla_a(fb) = f\nabla_a b + \#a(f)b, \qquad a, b \in \Gamma(E), \quad f \in C^\infty(B),
\] 
that is metric and torsion-free \cite[p.62]{boucetta2011riemannian}. 
It is characterised by Koszul's formula:
\begin{align*}
    2g(\nabla_a b, c) &= \#a(g(b, c)) + \#b(g(a,c)) - \#c(g(a,b))
    \\
    &\quad + g([a, b], c)  - g([b, c], a) + g([c, a], b),
\end{align*}
for \(a, b, c \in \Gamma(E)\). 
Its curvature and Ricci tensor are given by 
\[
R(a, b)c = \nabla_a \nabla_b c - \nabla_b \nabla_a c - \nabla_{[a, b]} c, 
\qquad
\ric_g(a, b) = \tr(c \mapsto R(c, a)b).
\]
Write \(\bar X = \sigma X\) and \(U = \iota U_{\rm ker}\) for \(X\in \mathfrak{X}(B)\) and \(U_{\rm ker} \in \Gamma(\ker \#)\). 
Denote the induced metric connection on \(\ker \#\) also by \(\nabla\). 
As in subsection \ref{subsec:Riemann-Geometry}, 
\[
\nabla_X U_{\rm ker} := \omega\nabla_{\bar X} U, \qquad \hat{\nabla}_{U_{\rm ker}} V_{\rm ker} := \omega \nabla_{U} V, \qquad L_X U_{\rm ker} := \omega \nabla_U \bar X.  
\]
Torsion freeness gives \(\nabla^{\rm ker} = \nabla  - L\), and O'Neill's integrability tensor is again \(A_{X, Y}= - \frac{1}{2}F_{X, Y}\).
\newline
\indent
The following proposition follows directly from Crainic--Fernandes \cite{crainic2003integbility} and Lie's third theorem. 
The correspondence of metrics, Levi-Civita connections, and curvature follows from Koszul's formula.

\begin{prop}[Local principal bundle realisation]\label{prop:local-principal-bundle-realise}
    Let \((E, B, [\cdot, \cdot], \#)\) be a transitive Lie algebroid equipped with a vector bundle metric \(g\). 
    Let \(U \subseteq B\) be open and contractible. 
    There is a simply-connected Lie group \(\mathsf{G}\) and an isomorphism of Lie algebroids over \(U\),
    \[
    E_U \cong T(U \times G)/\mathsf{G}
    \]
    under which \(g|_U\) corresponds to an invariant Riemannian metric on the principal bundle \(U \times \mathsf{G} \to U\). 
    The induced Levi-Civita connection, curvature and Ricci tensor agree with their definitions on \(E_U\), that is, they do not depend on the choice of realisation. 
\end{prop}

\begin{definition}[Deformations of Riemannian Lie algebroids (cf.~{\cite[Def.~1]{crainic2008deform}})]\label{def:deform-Riem-Lie-alg}
    Fix a vector bundle \(E \to B\). 
    A \textnormal{family of Riemannian Lie algebroids} over an interval \(I \subseteq \RR\) is a collection of Riemannian Lie algebroids 
    \[
    (E, B, [\cdot, \cdot]_t, \#_t, g_t)
    \]
    varying smoothly with \(t \in I\).      
    A family of isometric isomorphisms of Riemannian Lie algebroids
    \[
    \Phi(t) : (E, B, [\cdot, \cdot]_t, \#_t, g_t)
    \to (E', B', [\cdot, \cdot]_t', \#_t', g_t')
    \]
    varying smoothly with \(t\) is called an \textnormal{equivalence}. 
    If an equivalence between two families of Riemannian Lie algebroids exists, then they are \textnormal{equivalent}. 
    A \textnormal{deformation of a Riemannian Lie algebroid} \((E, B, [\cdot, \cdot], \#, g)\) is a family of Riemannian Lie algebroids over an interval \(I\) containing \(0\) with initial value
    \[
    (E, B, [\cdot, \cdot]_0, \#_0, g_0) = (E, B, [\cdot, \cdot], \#, g).
    \]
\end{definition}

\begin{example}\;
    \begin{enumerate}[(i)]
        \item A Ricci flow \(g(t)\) on a smooth manifold \(M\) is a family of Riemannian Lie algebroids with fixed bracket and fixed anchor:
        \[
        (TM, M, [\cdot,\cdot], \operatorname{Id}_{TM}, g(t)) 
        \]
        \item 
        The Ricci flow of left-invariant metrics on a Lie group with initial metric \(g_0\) is a deformation of a Riemannian Lie algebroid with fixed bracket and fixed anchor
        \[
        (\mathfrak{g}, *, [\cdot, \cdot], *, \langle \cdot, \cdot \rangle_t).
        \]
        Lauret's bracket flow is another deformation of the same Riemannian Lie algebroid but instead with fixed metric and fixed anchor:
        \[
        (\mathfrak{g}, *, \mu(t), *, \langle \cdot, \cdot \rangle_0).
        \]
        These two families of Riemannian Lie algebroids are equivalent \cite[Thm. 3.3]{lauret2013RicciFlowOf}. 
        The equivalence is constructed by solving either ODE, 
        \[
        \frac{d\Phi}{dt} = -\Phi \Ric_g, 
        \qquad 
        \frac{d\Phi}{dt} =  - \Ric_\mu \Phi,
        \]
        with initial value \(\Phi(0) = \operatorname{Id}\). 
        \item 
        The Ricci flow of invariant metrics on a twisted principal \(\mathsf{G}\)-bundle \(M \to B\) with initial metric \(g_0\) is a deformation of a Riemannian Lie algebroid with fixed bracket and fixed anchor:
        \[
        (TM/\mathsf{G}, B, [\cdot, \cdot], \pi_*, g(t)). 
        \]
    \end{enumerate}
\end{example}

\section{Bracket flows in higher cohomogeneity}

In this section, we prove \autoref{thm:B}. 
First, we write the Ricci flow of invariant metrics as a system of partial differential equations for the reduced data triple \((h, \omega, \check{g})\) \eqref{equ:reduce-data-bij}. 
We define the bracket flow on twisted principal bundles and prove its equivalence to the Ricci flow as deformations of Lie algebroids. 
In the final subsection, we discuss moment maps on the adjoint bundle, in preparation for the proof of \autoref{thm:A}.

\subsection{Ricci flow equations}
Let \(M \to B\) be a twisted principal \(\mathsf{G}\)-bundle. 
Equip \(M\) with an invariant Riemannian metric \(g\). 
Define the Ricci tensor and the Ricci endomorphism, 
\[
\ric_g:TM/\mathsf{G} \times TM/\mathsf{G} \to \RR, 
\qquad
    \Ric_g:TM/\mathsf{G} \to TM/\mathsf{G}. 
\]
With respect to the splitting given by the orthogonal complement, write 
\begin{align}
    \ric_g^{\rm ad}&: \adM \times \adM \to \RR, 
&
\ric_g^{\rm ad}(U_{\rm ad}, U_{\rm ad}) &= \ric_g( U, U),
\label{equ:Ricci-tensor-vert-decomp}
\\
\ric_g^{TB}&:TB \times TB \to \RR, 
&\ric_g^{TB}(X, X) &= \ric_g(\bar{X}, \bar{X}).
\label{equ:Ricci-tensor-hori-decomp}
\end{align}
Let \(\mathcal{V}, \mathcal{H}:TM/\mathsf{G} \to TM/\mathsf{G}\) be the projections of invariant vector fields into their vertical and horizontal components, as given in subsection \ref{subsec:Riemann-Geometry}. 
Decompose the Ricci endomorphism as the following endomorphisms of \(TM/\mathsf{G}\):
\begin{equation}\label{equ:Ricci-endo-decomp}
    \Ric_g^{\rm ad} =  \mathcal{V}\Ric_g \mathcal{V},
\qquad
\Ric_g^{\mathrm{ad}, TB} = \mathcal{V} \Ric_g \mathcal{H},
\qquad
\Ric_g^{TB, \mathrm{ad}} = \mathcal{H} \Ric_g \mathcal{V},
\qquad
\Ric_g^{TB} =  \mathcal{H}\Ric_g\mathcal{H}.
\end{equation}

\begin{prop}\label{prop:Ricci-flow-system}
On twisted principal bundles, the Ricci flow of invariant metrics \eqref{equ:Ricci-flow} is equivalent to the following system of partial differential equations for the triple \((h, \omega, \check{g})\):   
 \begin{align}
        \frac{\partial h}{\partial t} &= -2 \ric^{\rm ad}_g,
        \label{equ:Ricci-flow-h}
        \\
        \frac{\partial\omega}{\partial t} &= -2\omega \Ric_g^{\mathrm{ad}, TB},
        \label{equ:Ricci-flow-omega}
        \\
        \frac{\partial \check{g}}{\partial t} &= -2\ric^{TB}_g.
        \label{equ:Ricci-flow-g^}
\end{align}
where \(h\) is a vector bundle metric on \(\adM\), \(\omega\) is a principal connection, and \(\check{g}\) is a Riemannian metric on the base \(B\) \eqref{equ:reduce-data-bij}.
\end{prop}

\begin{proof}
    Equations \eqref{equ:Ricci-flow-h} and \eqref{equ:Ricci-flow-g^} are the restrictions of the Ricci flow equation to vertical and horizontal vectors. 
    For equation \eqref{equ:Ricci-flow-omega}, let \(g(t)\) be a Ricci flow of invariant metrics. 
    Take a vector field \(X\) on \(B\) and take its horizontal lift \(\bar{X}(t)\). 
    For all fixed vertical vectors \(U\), we know that \(g(t)(U,\bar{X}(t)) = 0\). 
    Differentiating this we get 
    \[
    0= -2\ric_g(U, \bar X) + g(U, \partial_t\bar X). 
    \]
    By definition \(\omega(t)(\bar X(t)) = 0\). Differentiating this we get 
    \[
    0 =  \frac{\partial\omega}{\partial t}(\bar X) + \omega(\partial_t \bar X) = \frac{\partial\omega}{\partial t}(\bar X) + 2\omega(\Ric_g \bar X),
    \]
    while \(\omega\) restricted to vertical vector fields remains unchanged under the Ricci flow. 
    This demonstrates equation \eqref{equ:Ricci-flow-omega}.
    Conversely, differentiating 
    \[
    g = \omega^*h + \pi^*\check{g} = h(\omega\cdot, \omega\cdot) + \check{g}(\pi_* \cdot, \pi_* \cdot)
    \]
    recovers the Ricci flow equation.
\end{proof}

\subsection{Bracket flow equations}
Let \(M \to B\) be a twisted principal \(\mathsf{G}\)-bundle. 
As in \ref{subsec:Lauret's}, given a Ricci flow of invariant metrics \(g(t)\) on \(M\), we apply Uhlenbeck's trick in the vertical directions by choosing a moving local frame \(\{U_j(t)\}\) of the adjoint bundle that solves
\[
 \partial_t U_j = \Ric_g^{\rm ad}(U_j).
\]
If \(\{U_i(t)\}\) is initially orthonormal, then it remains orthonormal along the Ricci flow. 
Indeed, the metric coefficients of the vertical directions are constant,
\begin{align*}
        \frac{\partial}{\partial t} g(U_i, U_j)
        &=\frac{\partial g}{\partial t} (U_i, U_j) + g(\partial_t U_i, U_j) + g(U_i, \partial_t U_j)
        \\
        &= -2\ric_g(U_i, U_j)
        + g(\Ric_g^{\mathrm{ad}}(U_i), U_j)
        + g( U_i, \Ric_g^{\rm ad}(U_j) )
        \\
        &= 0. 
    \end{align*}
At \(t = 0\), let \(\bar X_i\) be a local orthonormal horizontal frame.
Write 
\begin{align*}
    [U_i + \bar X_i, U_j +\bar X_j] &= [U_i, U_j] + \nabla_{ X_i}^{\rm ad} U_j - \nabla^{\rm ad}_{ X_j}U_i +[\bar X_i, \bar X_j]
    \\
    &= c_{ij}^k(X_i, X_j) U_k + \overline{[X_i, X_j]},
\end{align*}
where \(c_{ij}^k\) are the vertical structure coefficients of the Lie algebroid bracket in the moving frame \(\{\bar X_i, U_i(t)\}\). 
They are given by the following sum of structure coefficients of the adjoint bundle \(\mu_{ij}^k\), connection coefficients \(\omega_i^k \in \Omega^1(B)\), and curvature coefficients \(F^k \in \Omega^2(B)\):
\[
c_{ij}^k(X_i, X_j)= \mu_{ij}^k +\omega_j^k(X_i) - \omega_i^k(X_j) - F^k(X_i,X_j). 
\]
The time derivative of the vertical structure coefficients in this choice of moving frame is 
\begin{align*}
    \frac{\partial c_{ij}^k}{\partial t}(X_i, X_j) U_k &= [\Ric_g^{\rm ad} U_i, U_j + \bar X_j] + [U_i + \bar X_i, \Ric_g^{\rm ad} U_j] - \Ric_g^{\rm ad}[U_i + \bar X_i, U_j +\bar X_j]
    \\
    &= \delta_{[\cdot, \cdot]}(\Ric_g^{\rm ad})(U_i + \bar X_i, U_j +\bar X_j). 
\end{align*}


\begin{definition}[Bracket Flow]\label{def:bracket-flow}
    Let \(\pi:M \to B\) be a twisted principal \(\mathsf{G}\)-bundle equipped with an invariant Riemannian metric \(g_0\).
    Write 
    \(
    g_0 = 
    \omega_0^*h_0 + \pi^*\check{g}_0
    \)
    where \(h_0\) is a vector bundle metric on the adjoint bundle \(\adM \to B\), \(\omega_0\) is a principal connection, and \(\check{g}_0\) is a Riemannian metric on the base \(B\) (see \eqref{equ:reduce-data-bij}, \eqref{equ:Atiyah-metric-bij}). 
    A \textit{bracket flow} is a deformation 
    \begin{equation}\label{equ:bf:deform-algebroid}
    (TM/\mathsf{G}, B, [\cdot, \cdot]_{\mu}(t), \pi_*, g_{\mu}(t)), 
    \end{equation}  
    of the Riemannian Lie algebroid \((TM/\mathsf{G}, B, [\cdot, \cdot], \pi_*, g_0)\) with fixed anchor \(\pi_*\), such that the metric
    \[
    g_{\mu}(t) = \omega_{\mu}(t)^*h_0 + \pi^*\check{g}_{\mu}(t)
    \]
    and the bracket \([\cdot, \cdot]_{\mu}(t)\) solve the following system of partial differential equations 
    \begin{align}
     \pdv{[\cdot, \cdot]_\mu}{t} &= \delta_{[\cdot, \cdot]_\mu} ({\Ric}_{g_\mu}^{\rm ad} )
     \label{equ:bracket-flow-mu}
     \\
    \pdv{\omega_\mu}{t} &= -2\omega_\mu\Ric_{g_\mu}^{{\rm ad}, TB}  
    \label{equ:bracket-flow-omega}
    \\
    \pdv{\check{g}_\mu}{t}  &=  -2\ric_{g_\mu}^{TB}. 
    \label{equ:bracket-flow-g^}
    \end{align}
    Here, \(\ric_{g_{\mu}(t)}\) and \(\Ric_{g_{\mu}(t)}\) are the Ricci tensor and Ricci endomorphism of the evolving Riemannian Lie algebroids \((TM/\mathsf{G}, B, [\cdot, \cdot]_{\mu}(t), \pi_*, g_{\mu}(t))\), and the superscripts, \(\mathrm{ad}\) and \(TB\), denote the projections with respect to splittings by orthogonal complements of \(g_{\mu}(t)\) (\ref{equ:Ricci-tensor-vert-decomp} 
    --\ref{equ:Ricci-endo-decomp}). 
\end{definition}

\begin{remark}
Given a bracket flow on a twisted principal bundle \(M\to B\), write the restriction of the bracket to vertical vector fields as \(\mu(t)\). 
Indeed, \(\mu(t)\) is a section of \(\Lambda^2 \adM^*\otimes \adM\) that satisfies the Jacobi identity, so it is a Lie bracket on \(\adM\). 
Its evolution under the bracket flow is 
\[
\frac{\partial\mu}{\partial t} = \delta_\mu(\omega_\mu \Ric_{g_\mu}^{\rm ad}\iota). 
\]
The Atiyah sequence becomes the short exact sequence of Riemannian Lie algebroids, 
\[
0 \to (\adM, \mu(t), h_0)
\xrightarrow{\;\iota\;}
(TM/\mathsf{G}, [\cdot, \cdot]_{\mu}(t),  g_{\mu}(t))
\xrightarrow{\pi_*}
(TB, [\cdot,\cdot]_{TB}, \check{g}_{\mu}(t))
\to 0,   
\]
(with anchors to be understood as \(0, \pi_*\) and \( \operatorname{Id}_{TB}\)). 
Here, \((\adM, \mu(t), h_0)\) has a moving Lie bracket, but a fixed metric, 
while \((TB, [\cdot,\cdot]_{TB}, \check{g}_{\mu}(t))\) has a fixed Lie bracket, but a moving metric. 
\end{remark}

If the base were a point, then we recover Lauret's bracket flow for Lie groups,
\[
\frac{d\mu}{dt} = \delta_\mu(\Ric_{g_\mu}). 
\]
Now, we prove that the bracket flow is equivalent to the Ricci flow as deformations of Riemannian Lie algebroids in the sense of \autoref{def:deform-Riem-Lie-alg}.

\begin{theorem}[\autoref{thm:B}]
    \label{thm:precise-B}
    Let \(\pi:M \to B\) be a twisted principal \(\mathsf{G}\)-bundle equipped with an invariant metric \(g_0 = (h_0, \omega_0, \check{g}_0)\).
    Ricci flows of invariant metrics on \(M\) with initial metric \(g_0\) are in one-to-one correspondence with bracket flows (\autoref{def:bracket-flow}) with the same initial data. 
    There exists a unique smooth family of vector bundle automorphisms \(\Phi(t)\) of \(TM/\mathsf{G}\) over the identity that preserve the anchor and satisfy 
    \[
    g(t) = \Phi(t)^*g_{\mu(t)}. 
    \]
    The vector bundle automorphisms solve the fibrewise ODEs 
    \begin{enumerate}[(i)]
        \item \(\frac{\partial}{\partial t}\Phi = - \Phi\Ric_g^{\rm ad}, \quad \Phi(0) =\operatorname{Id}\),
        \item \(\frac{\partial}{\partial t}\Phi = - \Ric_{g_\mu}^{\rm ad}\Phi, \quad \Phi(0) =\operatorname{Id}\).
    \end{enumerate}
    Furthermore, the following conditions hold
    \begin{enumerate}[(i)]
    \setcounter{enumi}{2}
        \item 
        \(h = h_0(\Phi^{\rm ad}\cdot, \Phi^{\rm ad}\cdot)\)
        \item \([\cdot, \cdot]_\mu = \Phi[\Phi^{-1}\cdot, \Phi^{-1}\cdot]\) 
        \item \(\omega_\mu = \Phi^{\rm ad}\circ \omega \circ \Phi^{-1}\)
        \item \(\check{g}_\mu = \check{g}\)
    \end{enumerate}
    where \(\Phi \iota = \iota \Phi^{\rm ad}\) and \(\iota : \adM \to TM/\mathsf{G}\) is the inclusion. 
\end{theorem}

\begin{proof}
    Let \(g(t) = (h, \omega, \check{g})(t)\) be a Ricci flow. 
    Then by \autoref{prop:Ricci-flow-system},  \((h, \omega, \check{g})(t)\) solves equations (\ref{equ:Ricci-flow-h}--\ref{equ:Ricci-flow-g^}).
    Solve the linear ODE (i) fibrewise. 
    The solution \(\Phi(t)\) is smooth, unique, and has the same maximal time of existence as the Ricci flow. 
    Furthermore, \(\Phi(t): TM/\mathsf{G} \to TM/\mathsf{G}\) remains invertible as
    \[
    \pdv{}{t} \det\Phi = -(\det \Phi)\tr(\Ric_g^{\rm ad}), \qquad \det \Phi(0) = \det \operatorname{Id} = 1. 
    \]
    Since \(\pi_*\Ric_g^{\rm ad} = 0\) and \(\Phi(0) = \operatorname{Id}\), we find that \(\pi_*\Phi = \pi_*\), and so \(\Phi\) preserves the vertical distribution defining \(\Phi \iota = \iota \Phi^{\rm ad}\).
    The endomorphism \(\Phi^{\rm ad}(t) : \adM \to \adM\) satisfies 
    \[
    \frac{\partial \Phi^{\rm ad}}{\partial t} = -\Phi^{\rm ad} \omega \Ric_g \iota. 
    \]
    Consider the vector bundle metric \(h_0(\Phi^{\rm ad}\cdot, \Phi^{\rm ad}\cdot)\) on the adjoint bundle \(\adM \to B\). 
    Its time derivative is
    \begin{align*}
         \frac{\partial }{\partial t} h_0(\Phi^{\rm ad}\cdot, \Phi^{\rm ad}\cdot) &= h_0(\partial_t \Phi^{\rm ad}\cdot, \Phi^{\rm ad}\cdot) +h_0(\Phi^{\rm ad}\cdot, \partial_t \Phi^{\rm ad}\cdot) 
         \\
         &= - h_0(\Phi^{\rm ad} \omega \Ric_g \iota\cdot, \Phi^{\rm ad}\cdot) - h_0(\Phi^{\rm ad}\cdot, \Phi^{\rm ad} \omega \Ric_g \iota\cdot). 
    \end{align*}
    On the other hand, the vector bundle metric \(h(t)\) on \(\adM\) under Ricci flow solves \eqref{equ:Ricci-flow-h}. 
    Hence, 
    \begin{align*}
         \frac{\partial h}{\partial t} 
         &= -2h(\omega\Ric_g \iota \cdot, \cdot)
         \\
         &= -h(\omega\Ric_g \iota \cdot, \cdot) - h(\cdot,\omega\Ric_g \iota \cdot).
    \end{align*}
    By uniqueness, \(h = h_0(\Phi^{\rm ad}\cdot, \Phi^{\rm ad}\cdot)\), proving (iii). 
    Next, we show that if \(g(t)= (h, \omega, \check{g})(t)\) is a Ricci flow, then 
    \begin{equation}\label{equ:bf-pf-triple}
        \Phi[\Phi^{-1}\cdot, \Phi^{-1}\cdot], 
    \qquad  
    g_{\Phi} := (h_0, 
    \Phi^{\rm ad}\circ \omega \circ \Phi^{-1}, 
     \check{g})
    \end{equation}
    is a bracket flow. 
    The bracket \(\Phi[\Phi^{-1}\cdot, \Phi^{-1}\cdot]\) is compatible with the anchor \(\pi_*\Phi^{-1} = \pi_*\) and has initial value \([\cdot, \cdot]\).
    As \((\Phi^{\rm ad})^* h_0 = h\) and \(\pi_* \Phi^{-1} = \pi_*\), the vector bundle metric determined by the triple in \eqref{equ:bf-pf-triple} is isometric to the Ricci flow \(g(t)\):
    \[
    g_\Phi := (\Phi^{\rm ad}\circ \omega \circ \Phi^{-1})^*h_0 + \pi^*\check{g} = (\Phi^{-1})^*( \omega^*h + \pi^*\check{g}) = (\Phi^{-1})^*g(t). 
    \]
    Naturality of the Ricci curvature gives us 
    \[
    \Ric_{g_\Phi}  = \Phi\Ric_g \Phi^{-1}.
    \]
    The Ricci curvature blocks satisfy
    \[
    \Ric_{g_\Phi}^{\rm ad}  = \Phi\Ric_g^{\rm ad} \Phi^{-1}, 
    \qquad \Ric^{\mathrm{ad}, TB}_{g_\Phi}  = \Phi\Ric^{\mathrm{ad}, TB}_g \Phi^{-1}, 
    \qquad \ric^{TB}_{g_\Phi} = \ric^{TB}_g. 
    \]
    Then, the bracket evolves by
    \begin{align*}
        \frac{\partial}{\partial t}  \Phi[\Phi^{-1}\cdot, \Phi^{-1}\cdot] 
        &= -\Phi \Ric_g^{\rm ad}[\Phi^{-1}\cdot, \Phi^{-1}\cdot] 
        + \Phi[\Ric_g^{\rm ad}\Phi^{-1}\cdot, \Phi^{-1}\cdot] + \Phi[\Phi^{-1}\cdot, \Ric_g^{\rm ad}\Phi^{-1}\cdot]
        \\
        &= - \Ric_{g_\Phi}^{\rm ad}\Phi[\Phi^{-1}\cdot, \Phi^{-1}\cdot] 
        + \Phi[\Phi^{-1}\Ric_{g_\Phi}^{\rm ad}\cdot, \Phi^{-1}\cdot] + \Phi[\Phi^{-1}\cdot, \Phi^{-1}\Ric_{g_\Phi}^{\rm ad}\cdot]
        \\
        &=\delta_{\Phi[\Phi^{-1}\cdot, \Phi^{-1}\cdot]}( \Ric_{g_\Phi}^{\rm ad}),         
    \end{align*}
    the connection evolves by 
    \begin{align*}
        \frac{\partial}{\partial t}  (\Phi^{\rm ad}\circ \omega \circ \Phi^{-1} ) &= -(\Phi^{\rm ad}\Ric_g^{\rm ad})\circ \omega \circ \Phi^{-1}
        - \Phi^{\rm ad}\circ (2\omega\Ric_g^{\mathrm{ad},TB}) \circ \Phi^{-1} 
        +\Phi^{\rm ad}\circ \omega \circ (\Ric_g^{\rm ad} \Phi^{-1})
        \\
        &= - \Phi^{\rm ad}\circ (2\omega\Ric_g^{\mathrm{ad},TB}) \circ \Phi^{-1} 
        \\
        &= -2(\Phi^{\rm ad}\circ \omega \circ \Phi^{-1} ) \Ric_{g_\Phi}^{\mathrm{ad},TB},
    \end{align*}
    and the metric on the base evolves by
    \[
    \frac{\partial \check{g}}{\partial t} = -2\ric^{TB}_g = -2\ric^{TB}_{g_\Phi}. 
    \]
    Indeed, \eqref{equ:bf-pf-triple} is a bracket flow (\autoref{def:bracket-flow}).
    \newline
    \indent
    Conversely, let \((TM/\mathsf{G}, B, [\cdot, \cdot]_{\mu(t)}, \pi_*, g_{\mu(t)})\) be a bracket flow with initial metric \(g_0\).
    Solve the linear ODE (ii) fibrewise. 
    As before, the solution \(\Phi(t)\) is unique, smooth, invertible, preserves the anchor \(\pi_*\),  and has the same maximal time of existence as the bracket flow. 
    Observe that the pullback of the bracket is fixed:
    \[
    \frac{\partial}{\partial t} \Phi^{-1}[\Phi\cdot, \Phi\cdot]_\mu = \delta_{\Phi^{-1}[\Phi\cdot, \Phi\cdot]_\mu}(\Phi^{-1}\partial_t \Phi+ \Phi^{-1}\Ric_{g_\mu}^{\rm ad} \Phi) =0.
    \]
    Set \(g(t) = \Phi(t)^* g_\mu(t) 
    \), then 
    \[
    \Phi(t) : (TM/\mathsf{G}, [\cdot, \cdot], \pi_*, g(t))  \to (TM/\mathsf{G}, [\cdot, \cdot]_{\mu(t)}, \pi_*, g_{\mu(t)}) 
    \]
    is an isometric isomorphism of Riemannian Lie algebroids where \(g(t)\) is a Ricci flow:
    \begin{align*}
        \frac{\partial g}{\partial t} 
        &= \frac{\partial}{\partial t} \Phi^*g_\mu 
        \\
        &= \frac{\partial}{ \partial t} (h_0(\omega_\mu \Phi\cdot, \omega_\mu \Phi\cdot) +\check{g}_\mu(\pi_*\Phi\cdot, \pi_*\Phi \cdot))
        \\
        &= h_0((\partial_t \omega_\mu \Phi + \omega_\mu \partial_t \Phi )\cdot , \omega_\mu \Phi \cdot) 
        + h_0(\omega_\mu \Phi \cdot, (\partial_t \omega_\mu \Phi + \omega_\mu \partial_t \Phi )\cdot  ) +(\partial_t\check{g}_\mu)(\pi_*\cdot, \pi_* \cdot)
        \\
        &= -h_0( \omega_\mu (2 \Ric_{g_\mu}^{\mathrm{ad}, TB}+\Ric_{g_\mu}^{\rm ad}) \Phi \cdot , \omega_\mu \Phi \cdot) 
        -h_0(  \omega_\mu \Phi \cdot,\omega_\mu (2 \Ric_{g_\mu}^{\mathrm{ad}, TB}+\Ric_{g_\mu}^{\rm ad}) \Phi \cdot) 
        \\
        &\quad -2\ric_{g_\mu}^{TB}(\pi_*\cdot, \pi_*\cdot)
        \\
        &= -2h_0( \omega_\mu \Phi \Ric_g^{\mathrm{ad}, TB}\cdot, \omega_\mu \Phi\cdot)
        - h_0(\omega_\mu\Phi\Ric_g^{\rm ad} \cdot , \omega_\mu \Phi \cdot) 
        \\
        &\quad-2h_0(  \omega_\mu \Phi \cdot,\omega_\mu \Phi \Ric_g^{\mathrm{ad}, TB}\cdot)
        -h_0(\omega_\mu\Phi\cdot, \omega_\mu\Phi\Ric_g^{\rm ad} \cdot ) -2\ric_g^{TB}
        \\
        &= -2\ric_g^{TB\times \mathrm{ad}}  -2\ric_g^{\rm ad}- 2\ric_g^{ \mathrm{ad}\times TB} - 2\ric_g^{TB}
        \\
        &= -2\ric_g. 
    \end{align*}
    The second last equality follows from symmetry of the blocks. 
\end{proof}

\subsection{Moment maps on the adjoint bundle}
Let \(\mathsf{G}\) be a Lie group equipped with an invariant metric \(g\). 
Let \(\mu\) be the bracket of the Lie algebra \(\mathfrak{g}\), and let \(\langle \cdot, \cdot\rangle\) be the inner product on the Lie algebra given by evaluating \(g\) at the identity. 
The Ricci curvature endomorphism is 
\begin{equation}\label{equ:Ricci-homogeneous}
    \Ric_{g_\mu}^g = M_\mu^g - \frac{1}{2} B_\mu^g - S^g(\operatorname{ad}_\mu H_\mu^g).
\end{equation}
Here \(M_\mu^g\) is the self-adjoint endomorphism of \((\mathfrak{g}, \langle\cdot, \cdot \rangle)\) given by 
\[
\langle M_\mu^g U, U \rangle = - \frac{1}{2} \sum_{i,j} \langle \mu(U, V_i), V_j \rangle^2 + \frac{1}{4} \sum_{i, j} \langle \mu(V_i, V_j), U\rangle^2 
\]
where \(\{V_i\}\) is an orthonormal basis for \((\mathfrak{g}, \langle\cdot, \cdot \rangle)\).
Consider the action of \(\mathsf{GL}(\mathfrak{g})\) on \(\Lambda^2 \mathfrak{g}^* \otimes \mathfrak{g}\) given by \(\Phi \cdot \lambda = \Phi\lambda(\Phi^{-1}\cdot, \Phi^{-1}\cdot)\). 
Then, 
\[
4\langle M_\mu^g, \alpha \rangle = \frac{d}{dt}\Big|_{t = 0} \frac{1}{2}|\exp(t\alpha) \cdot \mu|^2 = -\langle \delta_\mu(\alpha), \mu\rangle
\]
for all \(\alpha \in {\rm End}(\mathfrak{g})\).
Because of this, we refer to this endomorphism \(M_\mu^g\) as the \textit{moment map} of the \(\mathsf{GL}(\mathfrak{g})\) action on \(\Lambda^2 \mathfrak{g}^* \otimes \mathfrak{g}\). 
The second term of \eqref{equ:Ricci-homogeneous} is the Killing form endomorphism given by 
\[
\langle B_\mu^g U, V\rangle = \tr \operatorname{ad}_\mu U \operatorname{ad}_\mu V.
\]
The third term of \eqref{equ:Ricci-homogeneous} is the symmetric part of the mean curvature operator \(\operatorname{ad}_\mu H_\mu^g\), where \(\langle H_\mu^g, U \rangle = \tr \operatorname{ad}_\mu U\) for all \(U \in \mathfrak{g}\). 
When \(\mathsf{G}\) is nilpotent, \(B_\mu^g\) and \(H_\mu^g\) vanish as all the eigenvalues of nilpotent matrices are zero, so \eqref{equ:Ricci-homogeneous} becomes just
\[
\Ric_{\mu}^g = M_\mu^g.
\]
\indent
Let \(M \to B\) be a principal \(\mathsf{G}\)-bundle equipped with an invariant metric \(g = (h, \omega, \check{g})\). 
Let \(\mu\) be the bracket on the adjoint bundle. 
On the fibres, \(\mathsf{G}\) acts freely and transitively, so the Ricci curvature of the fibres is given again by \eqref{equ:Ricci-homogeneous}.
With a fixed background metric \(h\) on \(\adM\), we write
\[
\hat{\Ric}_\mu = M_\mu - \frac{1}{2}B_\mu -S(\operatorname{ad}_\mu  H_\mu). 
\]
Here one can think of \(M_\mu \in {\rm End}(\adM)\) as the moment map of the gauge group action  (the vector bundle automorphisms over the identity) on sections of \(\Lambda^2 \adM^* \otimes \adM\):
\[
\Phi \cdot \lambda = \Phi \lambda(\Phi^{-1}\cdot, \Phi^{-1}\cdot), 
\]
where \(\Phi \in {\rm Aut}(\adM)\) and \(\lambda \in \Gamma(\Lambda^2 \adM^* \otimes \adM)\).
Indeed, 
\begin{equation}\label{equ:def-moment-map-implicit}
    h(M_\mu, \alpha) = -\frac{1}{4}h(\delta_\mu(\alpha), \mu),
\end{equation}
where \(\delta_\mu :{\rm End}(\adM) \to \Gamma(\Lambda^2 \adM^* \otimes \adM)\) is given by
\[
\delta_\mu(\alpha) =  \mu(\alpha \cdot, \cdot) + \mu( \cdot, \alpha \cdot) -\alpha\mu(\cdot, \cdot).
\]
More generally, let \([\cdot, \cdot]_\mu\) be a Lie algebroid bracket on \(TM/\mathsf{G}\) with anchor \(\pi_*\). 
Let  \(\mu \in \Gamma(\Lambda^2 \adM^* \otimes \adM)\) be the restriction of \([\cdot, \cdot]_\mu\).
A choice of principal connection on \(M\) determines an affine connection \(\nabla^{\rm ad\mu}\) on the adjoint bundle \(\adM\):
\[
\nabla^{\rm ad\mu}_X U_{\rm ad} = \omega[\bar X, U]_{\mu}. 
\]
By the Jacobi identity, we find that \(\mu\) is parallel with respect to this connection,   
\begin{align*}
    (\nabla^{\rm ad\mu}_X \mu)(U_{\rm ad}, V_{\rm ad}) &= \nabla^{\rm ad\mu}_X \mu(U_{\rm ad}, V_{\rm ad})  - \mu(\nabla^{\rm ad\mu}_X U_{\rm ad}, V_{\rm ad}) - \mu(U_{\rm ad}, \nabla^{\rm ad\mu}_X V_{\rm ad}) 
    \\
    &= [\bar{X}, [U, V]_\mu]_{\mu} - [[\bar{X}, U]_\mu, V]_{\mu} - [U, [\bar{X}, V]_\mu]_{\mu}
    \\
    &= 0.
\end{align*}
Equip the adjoint bundle \(\adM\) with a metric \(h\). 
If \(\nabla\) is a metric connection on \((\adM, h)\), then define \(L \in \Omega^1(B;\operatorname{End}(\adM))\) by  \(\nabla^{\rm ad\mu} = \nabla - L\).
As \(\nabla^{\operatorname{ad}\mu} \mu = 0\), we have
\[
 (\nabla_X \mu) = -\delta_\mu(L_X),
 \qquad
    (\nabla_{X, Y}^2 \mu)= -\delta_\mu((\nabla_X L)_Y) + \delta_{\delta(L_X)\mu}(L_Y).
\]


\begin{prop}\label{prop:1st-2nd-variation-moment-map}
    Let \(M \to B\) be a twisted principal bundle.
    Equip \(\adM \to B\) with a bundle metric \(h\). 
    Let \(\mu \in \Gamma(\Lambda^2 \adM^* \otimes \adM)\) satisfy the Jacobi identity. Then  
    \begin{align}
        h(\nabla_X M_\mu, \alpha) &= -\frac{1}{2}h(\delta_\mu(\alpha), \nabla_X \mu),
        \label{equ:moment-map-1st-derivative}
        \\
        h(\nabla_{X}\nabla_{Y} M_\mu, \alpha) &= -\frac{1}{2} h(\delta_{\nabla_X\mu}(\alpha), \nabla_Y \mu) - \frac{1}{2}h(\delta_\mu(\alpha), \nabla_{X} \nabla_{Y}\mu)
        \label{equ:moment-map-2nd-derivative}
    \end{align}
    for all \(h\)-symmetric \(\alpha \in \operatorname{End}(\adM)\) and for all vector fields \(X, Y\) on \(B\).
    Here \(M_\mu\) is defined by \eqref{equ:def-moment-map-implicit}.
\end{prop}

\begin{proof}
If \(\alpha \in {\rm End}(\adM)\) is \(h\)-symmetric, then \(\delta(\alpha)\in {\rm End}(\Lambda^2\adM^* \otimes \adM)\) is too. 
As \(\nabla_X \delta_\mu(\alpha)  = \delta_\mu(\nabla_X \alpha)+ \delta_{\nabla_X \mu}(\alpha)\), 
we can write
\begin{align*}
    4h(\nabla_X M_\mu, \alpha) &= 4Xh(M_\mu, \alpha) - 4h(M_\mu, \nabla_X \alpha)
    \\
    &= -Xh(\delta_\mu(\alpha), \mu) + h(\delta_\mu(\nabla_X \alpha), \mu)
    \\
    &= -h(\delta_{\nabla_X \mu}(\alpha), \mu) - h(\delta_\mu(\alpha),\nabla_X \mu)
    \\
    &= -2h(\delta_\mu(\alpha), \nabla_X \mu).
\end{align*}
For equation \eqref{equ:moment-map-2nd-derivative}, we substitute in equation \eqref{equ:moment-map-1st-derivative}:
\begin{align*}
    2h(\nabla_X\nabla_{Y} M_\mu, \alpha) &= 2Xh(\nabla_Y M_\mu, \alpha) - 2h(\nabla_Y M_\mu, \nabla_X \alpha)
    \\
    &= -Xh(\delta_\mu(\alpha), \nabla_Y\mu) +h(\delta_\mu(\nabla_X \alpha), \nabla_Y \mu)
    \\
    &= -h(\delta_{\nabla_X \mu}(\alpha),  \nabla_Y\mu) -h(\delta_\mu(\alpha), \nabla_{X} \nabla_{Y}\mu).
\end{align*}
\end{proof}

\section{Type III estimates}

In this section, we prove \autoref{thm:A}. 
To do this, we use the bracket flow.
Under the bracket flow, we prove Riemannian curvature and diameter estimates using the maximum principle. 
By \autoref{thm:B}, these estimates hold for the Ricci flow.

\subsection{Setup}
Let \(M\to S^1\) be a twisted principal \(\mathsf{G}\)-bundle.  Consider the Ricci flow of invariant metrics on \(M\) in the following gauge,
\begin{equation}\label{equ:gauge-ricci-flow}
    \frac{\partial g}{\partial t} = -2\ric_g + \mathcal{L}_N g, \qquad g(0) = g_0
\end{equation}
where \(N = -(\tr L)\bar X\) is the mean curvature vector of the fibres. 
This choice simplifies computations as it eliminates the first order terms in our maximum principle arguments. 
The Ricci curvature of invariant metrics on twisted principal bundles with one-dimensional base is given in \autoref{prop:Ricci-curv-dim1}. 
Then, in this gauge \eqref{equ:gauge-ricci-flow}, the equivalent formulation of the bracket flow (\autoref{def:bracket-flow})
\[
(TM/\mathsf{G}, B, [\cdot, \cdot]_{\mu}(t), \pi_*, g_{\mu}(t))
\]
where 
\[
g_\mu(t) = \omega_\mu(t)^*h_0 + \pi^*\check{g}_\mu(t),
\]
is the following system 
\begin{align}
     \frac{\partial[\cdot, \cdot]_\mu}{\partial t} &= \delta_{[\cdot, \cdot]_\mu} ( \hat{\Ric}_{\mu} -\nabla L),
     \label{equ:bf-dim-1-mu}
     \\
     \frac{\partial\omega_\mu}{\partial t}(\bar{X}) &= 2\omega_\mu ({\hat\nabla}^*L),
     \label{equ:bf-dim-1-omega}
     \\
     \frac{\partial\check{g}_\mu}{\partial t} &=  2|L|^2  \check{g}_\mu. 
     \label{equ:bf-dim-1-g^}
\end{align}
Here \(\hat{\Ric}_{\mu}\) is the Ricci endomorphism of the fibres. 
We do not use a subscript \(g_\mu\) as the Ricci curvature of the fibres only depends on the fixed vector bundle metric \(h\) on \(\adM\), and the moving Lie bracket \(\mu(t)\) on \(\adM\). 
Here, \(\bar X\) is the horizontal lift of a unit vector field \(X\) on the base, \(L = L_X\) is the shape operator, \(\nabla L = (\nabla_X L)_X\) is its covariant derivative, \(\hat{\nabla}\) is the Levi-Civita connection of the fibres, and its adjoint applied to the shape operator is 
\[
\hat{\nabla}^* L = \hat{\nabla}^*(L_X) =   -\sum_j (\hat{\nabla}_{U_j} L_X)_{U_j}
\]
for local orthonormal frames \(\{U_i\}\) of \((\adM, h_0)\). 
It follows from the proof of \autoref{thm:B} that (\ref{equ:bf-dim-1-mu}--\ref{equ:bf-dim-1-g^}) is equivalent to \eqref{equ:gauge-ricci-flow}. 
By \autoref{prop:Riemann-curv-dim1}, the Riemannian curvature of \((TM/\mathsf{G}, B, [\cdot, \cdot]_{\mu(t)}, \pi_*, g_{\mu(t)})\) is bounded above by
\begin{equation}
    |\Rm_{g_\mu}| \leq C(|\mu|^2 + |L|^2 +|\nabla L +L^2|)
\end{equation}
for some constant \(C>0\).
Here \(\mu(t)\) is the Lie bracket on \(\adM\) given by the restriction of \([\cdot, \cdot]_\mu\) to vertical vector fields.
In the following subsections, we estimate each of these terms.
First, we compute the evolution of the shape operator, the Levi-Civita connection \(\nabla\) on \(\adM\), and
the connection \(\nabla^{\rm ad\mu}\) on \(\adM\) induced by the principal connection:
\[
\nabla^{\rm ad\mu}_X U_{\rm ad} = \omega_\mu[\bar{X}, U]_\mu, 
\]

\begin{lemma}\label{lemma:bf-connections-evolve}
    Let \(M \to S^1\) be a twisted principal \(\mathsf{G}\)-bundle. 
    Fix a vector bundle metric \(h\) on the adjoint bundle \(\adM\) and let \(([\cdot, \cdot]_\mu, \omega_\mu, \check{g}_\mu)(t)\) be a bracket flow on \(M\) solving (\ref{equ:bf-dim-1-mu}--\ref{equ:bf-dim-1-g^}).
    The connections and the shape operator evolve by 
    \begin{align*}
        \partial_t (\nabla^{\rm ad\mu}_X) &=\nabla^{\rm ad\mu}_X(\hat{\Ric}_{\mu} -\nabla L)-2\operatorname{ad}_\mu ({\hat\nabla}^*L) -|L|^2\nabla_{X}^{\rm ad
        \mu},
        \\
        \partial_t L&= -\nabla_X(\hat{\Ric}_{\mu} -\nabla L)+2\operatorname{ad}_\mu ({\hat\nabla}^*L)^{\rm sym} -|L|^2 L,
        \\
        \partial_t(\nabla_X) &=-[L, \hat{\Ric}_{\mu} -\nabla L]-2\operatorname{ad}_\mu ({\hat\nabla}^*L)^{\rm skew} -|L|^2\nabla_{X}.
    \end{align*}
    Here \(X\) is a unit vector field on \(S^1\), \([\cdot, \cdot]\) is the commutation of endomorphisms of \(\adM\), and the superscripts sym and skew denote the symmetric and skew-symmetric components of an endomorphism of \(\adM\) with respect to \(h\).   
\end{lemma}

\begin{proof}
    If \(\bar{X}\) is horizontal and \(X\) is unit at \(t = 0\), for this to persist, the vector field \(\bar{X}\) must evolve by 
    \[
    \partial_t \bar{X} = -2\hat\nabla^* L
    -|L|^2\bar{X}, 
    \qquad
    \partial_t X = - |L|^2 X
    \]  
    Fix a time-independent vertical vector \(U\). 
    To take the time derivative of \(\nabla^{\rm ad\mu}_X U_{\rm ad}\), 
    we write 
    \[
    \nabla^{\rm ad\mu}_X U_{\rm ad} = \omega_\mu[\bar{X}, U]_\mu, 
    \]
    so then we get 
    \begin{align*}
         \partial_t \nabla^{\rm ad\mu}_X U_{\rm ad} &= \frac{\partial \omega_\mu}{\partial t} [\bar{X}, U]_\mu
         +
         \omega_\mu \frac{\partial [\cdot, \cdot]_\mu}{\partial t}(\bar X, U) 
         +
         \omega_\mu[\partial_t \bar{X}, U]_\mu
         \\
         &=  -2\mu({\hat\nabla}^*L, U_{\rm ad}) +\left(\nabla^{\rm ad\mu}_X(\hat{\Ric}_{\mu} -\nabla L)\right) U_{\rm ad}-|L|^2\nabla_{X}^{\rm ad\mu} U_{\rm ad}. 
    \end{align*}
    With respect to a chosen \(h\)-orthonormal frame, the connection matrix of \(\nabla^{\rm ad}\) splits into a skew-symmetric part and a symmetric part. 
    The skew-symmetric part is the connection matrix of the metric connection \(\nabla\). 
    The symmetric part is \(-L\). 
    Because \(h\) does not depend on time, we get 
    \[
    \partial_t \nabla_X +|L|^2\nabla_X = (\partial_t \nabla^{\rm ad\mu}_X+|L|^2\nabla_X^{\rm ad\mu} )^{\rm skew}, 
    \qquad
    \partial_t L + |L|^2L= -(\partial_t \nabla^{\rm ad\mu}_X+ |L|^2\nabla_X^{\rm ad\mu})^{\rm sym} 
    \]
    in a fixed orthonormal frame. 
    In particular, we get the splitting into symmetric and skew-symmetric parts
    \[
    \nabla^{\rm ad\mu}_X(\hat{\Ric}_{\mu} -\nabla L) = \nabla_X(\hat{\Ric}_{\mu} -\nabla L) -[L, \hat{\Ric}_{\mu} -\nabla L]. 
    \]
\end{proof}

\begin{lemma}
    Let \(M \to S^1\) be a twisted principal \(\mathsf{G}\)-bundle. 
    Fix a vector bundle metric \(h\) on the adjoint bundle \(\adM\),  and let \(([\cdot, \cdot]_\mu, \omega_\mu, \check{g}_\mu)(t)\) be a bracket flow on \(M\) solving  (\ref{equ:bf-dim-1-mu}--\ref{equ:bf-dim-1-g^}). 
    Take a vector \((U_{\rm ad})_x \in \adM_x\) and parallel transport it to a local section \(U_{\rm ad}\) via the metric connection \(\nabla(t)\) induced by the Levi-Civita connection.
    If \(U_{\rm ad}(t)\) evolves by
    \[
    \partial_t U_{\rm ad}= Q U_{\rm ad}
    \]
    where \(Q(t)\) is a skew-symmetric endomorphism of \((\adM, h)\) solving 
    \[
    \nabla_X Q = [L, \hat{\Ric}_{\mu} - \nabla L] + 2(\operatorname{ad}_\mu ({\hat{\nabla}}^* L))^{\rm skew}
    , \qquad Q(x, t) = 0, 
    \]
    then \(U_{\rm ad}(t)\) is \(\nabla(t)\)-parallel.
\end{lemma}

\begin{proof}
    From \autoref{lemma:bf-connections-evolve},
    \begin{align*}
        \partial_t  (\nabla_X U_{\rm ad}) &= \partial_t (\nabla_X)U_{\rm ad} + \nabla_X \partial_t U_{\rm ad}
        \\
        &=\partial_t (\nabla_X)U_{\rm ad} +\nabla_X (Q U_{\rm ad})
        \\
        &= \partial_t (\nabla_X)U_{\rm ad} +(\nabla_X Q) U_{\rm ad} +Q \nabla_X  U_{\rm ad}
        \\
        &=   -|L|^2\nabla_{X} U_{\rm ad} +Q \nabla_X  U_{\rm ad}.
    \end{align*}
    This is a homogeneous linear ODE for \(\nabla_X U_{\rm ad}\) with initial value  \(\nabla_X U_{\rm ad}(0) = 0\), and thus, it will be zero for all time. 
    Indeed, if \(U_{\rm ad}(t)\) is initially parallel, it will remain parallel along the bracket flow. 
\end{proof}

\subsection{Lie bracket estimate}

\begin{prop}\label{prop:bracket-bound}
    Let \(\mathsf{N}^n\) be a nilpotent Lie group. 
    Let \(M \to S^1\) be a twisted principal \(\mathsf{N}\)-bundle.
    Fix a vector bundle metric \(h\) on \(\adM\).  
    If \(([\cdot, \cdot]_\mu, \omega, \check{g})(t)\) is a bracket flow on \(M\) solving  (\ref{equ:bf-dim-1-mu}--\ref{equ:bf-dim-1-g^}), then for \(t>0\), the Lie bracket on the adjoint bundle \(\adM\) satisfies
    \[
    |\mu |^2 \leq \frac{2n}{t}.
    \]
\end{prop}

\begin{proof}
    Take \(X(t)\) to be a unit vector field on \((S^1, \check{g}(t))\).
    All norms are given with respect to \(h\). 
    The square norm of the Lie bracket \(|\mu|^2\) is a function of time and the base \(S^1\).
    By the Jacobi identity, the bracket \(\mu\) is \(\nabla^{\rm ad\mu}\)-parallel, 
    \[
    0=\nabla^{\rm ad\mu}_X \mu = \nabla_X \mu + \delta_\mu(L).  
    \]
    In a \(\nabla(t)\)-parallel frame of \(\adM\), the Laplacian is
    \begin{align*}
    \Delta_{S^1} |\mu|^2 &= X^2 |\mu|^2
    \\
    &= 2|\nabla_X \mu|^2 + 2h(\nabla_X \nabla_X \mu, \mu) 
    \\
    &= 2|\delta_\mu(L)|^2 + 2h(\delta_{\delta_\mu(L)}(L)-\delta_\mu\nabla L, \mu),
    \\
    &= 4|\delta_\mu(L)|^2 - 2h(\delta_\mu\nabla L, \mu).
    \end{align*}
    In a fixed frame, the time derivative is
    \begin{align*}
        \pdv{}{t} |\mu|^2  &= 2h(\partial_t \mu, \mu)
        \\
        &= 2h(\delta_\mu(\hat{\Ric}_\mu), \mu) -  2h(\delta_\mu\nabla L, \mu).
    \end{align*}
    Hence, the heat operator is
    \begin{align*}
        \left( \frac{\partial}{\partial t} -\Delta_{S^1} \right) |\mu|^2 &= 2h(\delta_\mu(\hat{\Ric}_\mu), \mu) - 4|\delta_\mu(L)|^2
        \\
        &= -8h(M_\mu, \hat{\Ric}_\mu) -4|\delta_\mu(L)|^2
    \end{align*}
    where \(M_\mu\) is the moment map.
    As \(\mathsf{N}\) is nilpotent, the Ricci curvature of the fibres is \(\hat{\Ric}_\mu = M_\mu\), while the scalar curvature of the fibres is \(\hat{R}_\mu = \langle M_\mu, \operatorname{Id}\rangle = -\frac{1}{4}|\mu|^2\). 
    By Cauchy--Schwarz, 
    \[
    \left( \frac{\partial}{\partial t} -\Delta_{S^1} \right) |\mu|^2
        \leq -8|M_\mu|^2  \leq
    -\frac{8}{n}\langle M_\mu, \operatorname{Id} \rangle^2
    = -\frac{1}{2n} |\mu|^4.
    \]
    By the maximum principle,
    \[
    |\mu|^2 \leq \frac{2n}{t}.
    \]
\end{proof}

\subsection{Shape operator estimate}

\begin{prop}\label{prop:shape-operator-estimate}
    Let \(\mathsf{N}\) be a nilpotent Lie group. 
    Let \(M \to S^1\) be a twisted principal \(\mathsf{N}\)-bundle.
    Fix a vector bundle metric \(h\) on \(\adM\).   
    If \(([\cdot, \cdot]_\mu, \omega, \check{g})(t)\) is a bracket flow on \(M\) solving (\ref{equ:bf-dim-1-mu}--\ref{equ:bf-dim-1-g^}), then for \(t> 0\),
    \[
    |L|^2 \leq \frac{1}{2t}. 
    \]
\end{prop}

\begin{proof}
    Take \(X(t)\) to be a unit vector field on \((S^1, \check{g}(t))\).
    All norms are given with respect to \(h\). 
    In a \(\nabla(t)\)-parallel frame, the Laplacian of the square norm is 
    \[
    \Delta_{S^1}|L|^2 =  X^2 |L|^2 = 2|\nabla L|^2 +2h(\nabla^2 L, L).
    \]
    From \autoref{lemma:bf-connections-evolve}, in a fixed frame, the time derivative is
    \begin{align*}
        \frac{\partial}{\partial t} |L|^2
        &= 2h(\partial_t L, L)
        \\
        &=-2h(\nabla\hat{\Ric}_\mu , L)
        +2h(\nabla^2 L , L) + 4 h(\operatorname{ad}_\mu({\hat\nabla}^*L), L) - 2|L|^4.
    \end{align*}
    If \(\hat{\nabla}\) is the Levi-Civita connection of the fibres, by torsion-freeness, 
    \[
    \operatorname{ad}_\mu({\hat\nabla}^*L) = \hat{\nabla}_{{\hat\nabla}^*L} - \hat{\nabla}_\cdot ({\hat\nabla}^*L).
    \]
    Since \(\hat{\nabla}_{{\hat\nabla}^*L} \in {\rm End}(\adM)\) is skew-symmetric, we get 
    \[
    h(\operatorname{ad}_\mu({\hat\nabla}^*L), L)= -h(\hat{\nabla}_\cdot {\hat\nabla}^*L, L) =  -|{\hat\nabla}^*L|^2 + h({\hat\nabla}^*L, LH),
    \]
    where \(H = \sum_i \hat\nabla_{U_i} U_i\) is the mean curvature vector of the fibres as homogeneous spaces.
    The heat operator is then 
    \begin{equation}\label{equ:heat-shape-operator}
        \left(\frac{\partial}{\partial t} - \Delta_{S^1}\right) |L|^2
        = -2h(\nabla\hat{\Ric}_\mu , L)
        -2|\nabla L|^2  -4|{\hat\nabla}^*L|^2 + 4h({\hat\nabla}^*L, LH)- 2|L|^4.
    \end{equation}
    As \(\mathsf{N}\) is nilpotent, \(H =0\) and \(\hat{\Ric}_\mu = M_\mu\). 
    Moreover, by \autoref{prop:1st-2nd-variation-moment-map}
    \[
    -2h(\nabla M_\mu , L) = -|\delta_\mu(L)|^2. 
    \]
    The desired inequality follows from the maximum principle.
\end{proof}

\subsection{Diameter estimate}

\begin{prop}
    Let \(\mathsf{N}\) be a nilpotent Lie group. 
    Let \(M \to S^1\) be a twisted principal \(\mathsf{N}\)-bundle.
    Fix a vector bundle metric \(h\) on \(\adM\).   
    If \(([\cdot, \cdot]_\mu, \omega, \check{g})(t)\) is a bracket flow on \(M\) solving (\ref{equ:bf-dim-1-mu}--\ref{equ:bf-dim-1-g^}), then for \(t> 0\)
    \[
    \diam(S^1, \check{g}(t)) \leq C_0\sqrt{1+t},
    \]
    where \(C_0\) is some positive constant depending only on the initial metric. 
\end{prop}

\begin{proof}
    By equation \eqref{equ:heat-shape-operator} in the proof of \autoref{prop:shape-operator-estimate},  
    \[
    |L|^2 \leq \frac{1}{2t+c_0}
    \]
    where \(c^{-1}_0 =\max\{1,  \sup_{S^1}|L(0)|^2\}\). 
    The derivative of the diameter is bounded by
    \[
    \frac{d}{dt} \diam(S^1, \check{g}(t)) = \frac{d}{dt} \frac{1}{2}\int_{S^1} ds= \frac{1}2{}\int_{S^1}|L|^2\, ds \leq \frac{1}{2t+c_0} \diam(S^1, \check{g}(t)). 
    \]
    By integrating, we get
    \[
    \diam(S^1, \check{g}(t)) \leq \diam(S^1, \check{g}(0)) \sqrt{1+ \frac{2t}{c_0}}
    \leq C_0\sqrt{1+t}.
    \]
\end{proof}

\subsection{Jacobi operator estimate}

\begin{prop}\label{prop:Jacobi-operator-estimate}
    Let \(\mathsf{N}^n\) be a nilpotent Lie group. 
    Let \(M \to S^1\) be a twisted principal \(\mathsf{N}\)-bundle.
    Fix a vector bundle metric \(h\) on \(\adM\).   
    If \(([\cdot, \cdot]_\mu, \omega, \check{g})(t)\) is a bracket flow on \(M\) solving (\ref{equ:bf-dim-1-mu}--\ref{equ:bf-dim-1-g^}), then for \(t> 0\)
    \[
    |L^2 +\nabla L| \leq \frac{C}{t}
    \]
    where \(C =C(n) > 0\) is some constant depending only on dimension.  
\end{prop}

\begin{proof}[Proof of \autoref{prop:Jacobi-operator-estimate}]
    We apply the maximum principle to \(|L|^2 + |\nabla L|\). 
    First, we compute the heat operator of \(|\nabla L|^2\). 
    Take \(X(t)\) to be a unit vector field on \((S^1, \check{g}(t))\).
    All norms are with respect to \(h\). 
    In a \(\nabla(t)\)-parallel frame, the Laplacian is 
    \[
    \Delta_{S^1}|\nabla L|^2 =  X^2 |\nabla L|^2 = 2|\nabla^2 L|^2 +2h(\nabla^3 L,\nabla L). 
    \]
    In a fixed frame, the time derivative is 
    \[
      \frac{\partial}{\partial t}|\nabla L|^2 = 2h(\partial_t \nabla L, \nabla L) = 2h((\partial_t \nabla_X)L, \nabla L)
    + 2h(\nabla_X \partial_t L, \nabla L).
    \]
    By \autoref{lemma:bf-connections-evolve}, the connection evolves by
    \[
    (\partial_t \nabla_X)L = - [[L,\hat{\Ric}_\mu- \nabla L], L] - 2[\operatorname{ad}_\mu ({\hat\nabla}^*L)^{\rm skew}, L] - |L|^2 \nabla L,
    \]
    while the shape operator evolves by
    \[
    \nabla_X \partial_t L =  -\nabla^2\hat{\Ric}_{\mu} + \nabla^3 L +2\nabla \operatorname{ad}_\mu ({\hat\nabla}^*L)^{\rm sym} -X(|L|^2) L -|L|^2 \nabla L.
    \]
    In a fixed frame, the time derivative is     
    \begin{align*}
        \frac{\partial}{\partial t}|\nabla L|^2 &= -2h([[L,\hat{\Ric}_\mu], L], \nabla L) -2|[L, \nabla L]|^2
        -4h([\operatorname{ad}_\mu ({\hat\nabla}^*L)^{\rm skew}, L] , \nabla L)
    \\
    &\quad-2h(\nabla^2\hat{\Ric}_{\mu}, \nabla L)
    +2h(\nabla^3 L, \nabla L)
    +4h(\nabla_X\operatorname{ad}_\mu ({\hat\nabla}^*L), \nabla L) 
    \\
    &\quad -4h(L, \nabla L)^2 -4|L|^2 |\nabla L|^2. 
    \end{align*}
    Since \(\mathsf{N}^n\) is nilpotent, \(\hat{\Ric}_\mu =M_\mu\), so \(|\hat{\Ric}_\mu| \leq  C|\mu|^2\)  for some \(C = C(n) >0\).  
    By Koszul's formula \(|{\hat\nabla}^*L| \leq C|\mu||L|\), and so 
    \[
    -2h([[L,\hat{\Ric}_\mu], L], \nabla L) -4h([\operatorname{ad}_\mu ({\hat\nabla}^*L)^{\rm skew}, L] , \nabla L)\leq C|\mu|^2|L|^2|\nabla L|.
    \]
    By \autoref{prop:1st-2nd-variation-moment-map}, the second spatial variation of the moment map is
    \[
    -2h(\nabla^2 M_\mu, \nabla L) 
        = h(\delta_{\delta_\mu(L)}(\nabla L), \delta_\mu(L)) 
    -|\delta_\mu\nabla L|^2
    + h(\delta_\mu(\nabla L), \delta_{\delta_\mu(L)}(L)). 
    \]
    Clearly, \(|\delta_\mu(L)| \leq 3|L||\mu|\),
    so by Cauchy--Schwarz
    \[
    -2h(\nabla^2 M_\mu, \nabla L)  \leq 54|\mu|^2|L|^2|\nabla L| -|\delta_\mu(\nabla L)|^2 . 
    \]
    Next, as \(\nabla \mu = -\delta_\mu(L)\), 
    \begin{align*}
        |\nabla {\hat\nabla}^*L| &\leq C(|\nabla \mu||L| +|\mu||\nabla L|) 
        \\
        &\leq C(3|\mu||L|^2 + |\mu||\nabla L|).
    \end{align*}
    We know that 
    \[
    \nabla_X \operatorname{ad}_\mu ({\hat\nabla}^*L)
    = \operatorname{ad}_{\nabla_X \mu}({\hat\nabla}^*L) +\operatorname{ad}_\mu(\nabla_X {\hat\nabla}^*L), 
    \]
    so
    \begin{align*}
        | \nabla_X \operatorname{ad}_\mu ({\hat\nabla}^*L)|&\leq |\nabla \mu||{\hat\nabla}^*L| +|\mu| |\nabla_X {\hat\nabla}^*L|
        \\
        &\leq C|\mu|^2(|L|^2 +|\nabla L|).
    \end{align*}
    By Propositions \ref{prop:bracket-bound} and \ref{prop:shape-operator-estimate}, the heat operator of the square norm is bounded by 
    \begin{align*}
        \left(\frac{\partial}{\partial t}-\Delta_{S^1}\right)|\nabla L|^2
        &\leq -|\delta_\mu(\nabla L)|^2 +C|\mu|^2|\nabla L|(|L|^2 +|\nabla L|) 
        -4h(L, \nabla L)^2 -4|L|^2 |\nabla L|^2-2|\nabla^2 L|^2
        \\
        &\leq -|\delta_\mu(\nabla L)|^2 -4h(L, \nabla L)^2-4|L|^2 |\nabla L|^2 -2|\nabla^2 L|^2
        + \frac{C}{t} |\nabla L|(|L|^2 +|\nabla L|).  
    \end{align*}
    The heat operator of the norm is then
    \[
    \left(\frac{\partial}{\partial t}-\Delta_{S^1}\right)|\nabla L|
    \leq 
    -\frac{|\delta_\mu(\nabla L)|^2}{2|\nabla L|} -\frac{2h(L, \nabla L)^2}{|\nabla L|} -2|L|^2 |\nabla L|
        + \frac{C}{t} (|L|^2 +|\nabla L|),
    \]
    and by adding \eqref{equ:heat-shape-operator} and positive terms, we get
    \[
      \left(\frac{\partial}{\partial t}-\Delta_{S^1}\right)(|L|^2+|\nabla L|)
      \leq -2|\nabla L|^2 -2|L|^2 |\nabla L| -2|L|^4 +\frac{C}{t}(|L|^2 +|\nabla L|).
    \]
    By the maximum principle, \(|L|^2+|\nabla L|\) is bounded by solutions to 
    \[
    \frac{dy}{dt} = -y^2 + \frac{C}{t} y,
    \]
    with initial value greater than  \(|L|^2+|\nabla L|\) at \(t = \varepsilon >0\).  
    Fix some \(c_0 > \sup_{S^1 \times [0, \varepsilon]}(|L|^2+|\nabla L|)\). 
    The solution to the ODE with initial value \(y(\varepsilon)= c_0\) is
    \[
    y = \frac{C+1}{t +\left(\frac{C+1}{c_0} - \varepsilon \right)(\frac{t}{\varepsilon})^{-C}}.
    \]
    For small \(\varepsilon>0\), we get \[
    y \leq \frac{C+1}{t}
    \]
    where \(C = C(n) >0\).
    The desired inequality for  \(|\nabla L + L^2|\) follows from the triangle inequality. 
\end{proof}

This completes the proof of \autoref{thm:A}.

\section{Convergence}

In this section, we prove \autoref{thm:C} and \autoref{coro:D}. 
First, we give compactness results for Ricci flows of invariant metrics on twisted principal bundles with simply-connected nilpotent fibres and \(S^1\) base.

\subsection{Compactness}

\begin{theorem}\label{thm:inv-RF-compactness}
    Fix \(n \in \mathbb{N}\) and \(K, d >0\). 
    Let \(I\subseteq \RR\) be an open interval containing \(0\). 
    Let \(M_k \to S^1\) be a sequence of twisted principal \(\mathsf{N}_k\)-bundles where \(\mathsf{N}_k\) is nilpotent, simply-connected and \(\dim \mathsf{N}_k = n\). 
    For every \(k \in \NN\), let \(g_k(t)\) be a normalised Ricci flow of invariant metrics on \(M_k\) defined for \(t \in I\) such that 
    \[
    |\Rm_{g_k(t)}|_{g_k(t)} \leq K, \qquad \diam(S^1, \check{g}_k(t)) \leq d. 
    \]
    Take a sequence of points in the universal covers \(\tilde{p}_k \in \tilde{M}_k\). 
\begin{enumerate}[(a)]
    \item 
    Then, the lifted flows \(\tilde{g}_k(t)\) subconverge to a normalised Ricci flow
    \[
    (\tilde{M}_k, \tilde{g}_k(t), \tilde{p}_k) \to (\tilde{M}_\infty, \tilde{g}_\infty(t), \tilde{p}_\infty)
    \]
    in the pointed smooth Cheeger--Gromov topology on compact time intervals in \(I\), and \(\tilde{M}_\infty\) is a principal \(\mathsf{N}_\infty\)-bundle over \(\RR\) for some simply-connected nilpotent Lie group \(\mathsf{N}_\infty\), and \(\tilde{g}_\infty(t)\) is \(\mathsf{N}_\infty\)-invariant. 
    \item 
    If in addition, there exists a constant \(c>0\) such that for all \(k \in \NN\),
    \[
    0 < c\leq \diam(S^1, \check{g}_k(t)),
    \]
    then \((\tilde{M}_\infty, \tilde{g}_\infty(t))\) is the lift to the universal cover of some normalised Ricci flow of invariant metrics \(g_\infty(t)\) on a twisted principal \(\mathsf{N}_\infty\)-bundle \(M_\infty \to S^1\). 
    \item 
    Assuming also that there exists a constant \(l >0\) and generators \(\tau_k\) of the deck transformation groups \(\tilde{M}_k \to M_k\) such that 
    \[
    \operatorname{dist}_{\tilde{g}_k(0)}(\tilde{p}_k, \tau_k(\tilde{p}_k)) \leq l,
    \]
    then, the flows \(g_k(t)\) subconverge
    \[
    (M_k, g_k(t), p_k) \to (M_\infty, g_\infty(t), p_\infty)
    \]
    in the pointed smooth Cheeger--Gromov topology on compact time intervals in \(I\). 
    Here \(M_\infty \to S^1\) is a twisted principal \(\mathsf{N}_\infty\)-bundle and \(g_\infty(t)\) is a normalised Ricci flow of invariant metrics whose lift to the universal cover is \((\tilde{M}_\infty, \tilde{g}_\infty(t))\). 
\end{enumerate}

\end{theorem}

\begin{proof}[Proof of \autoref{thm:inv-RF-compactness}] 
    The metrics \(g_k(t)\) are complete as the base is compact and the metrics on the fibres are homogeneous. 
    The lifts to the universal cover \(\tilde{g}_k(t)\) are complete too. 
    \newline
    \newline
    \indent
    \textbf{Part (a).}
    First, we prove an injectivity radius bound on \((\tilde{M}_k, \tilde{g}_k(0))\) and then apply Hamilton's compactness theorem \cite{hamilton1996compactness}.  
    For brevity, write \(\tilde{g}_k := \tilde{g}_k(0)\), and use horizontal transport to identify 
    \begin{equation}\label{equ:pf-compact-identify}
        (\tilde{M}_k, \tilde{g}_k, \tilde{p}_k) \cong (\RR \times \mathsf{N}_k, ds^2 +h_k(s), (0,\mathsf{e})),
    \end{equation}
    where \(h_k(s)\) is a one-parameter family of right-invariant metrics on \(\mathsf{N}_k\).
    The Lie exponential maps of \(\mathsf{N}_k\) are diffeomorphisms. 
    As \(\operatorname{ad}_X\) is nilpotent for all \(X \in \mathfrak{n}_k\), the differential after right translation back to the identity 
    \[
    (R_{\exp(-X)})_*d(\exp)_X 
    = \frac{e^{\operatorname{ad}_X}-\operatorname{Id}}{\operatorname{ad}_X}=\operatorname{Id} + \frac{1}{2}\operatorname{ad}_X + \frac{1}{6} (\operatorname{ad}_X)^2 + \cdots
    \]
    \cite[Theorem 2.14.3]{varadarajan1984lie} 
    is unipotent. 
    For all \(r>0\),
    \[
    \Vol_h\exp (B(0, r)) = \int_{B(0, r)} |\det d(\exp)_X | \;dx = \Vol_{(\mathfrak{n}, h_\mathsf{e})} B(0,r) = \omega_n r^n
    \]
    where \(dx\) is the Euclidean measure on \((\mathfrak{n}, h_\mathsf{e})\), and \(\omega_n\) is the volume of the \(n\)-dimensional Euclidean unit ball.
    Given \(X \in \mathfrak{n}\), by right-invariance the length of the curve \(\gamma_X(t) = \exp(tX)\) satisfies
    \[
    \operatorname{dist}_{h}(\mathsf{e}, \exp X)\leq \operatorname{Length}_h(\gamma_X) = \int_0^1 |\dot{\gamma}_X(t)|_h \; dt
    = \int_0^1 |(R_{\exp(tX)})_* X|\; dt = \int_0^1 |X|_{h_\mathsf{e}}dt = |X|_{h_\mathsf{e}}.  
    \]
    In particular, 
    \[
    \Vol_h B(\mathsf{e}, r) \geq \Vol_h\exp (B(0, r)) = \omega_n r^n .
    \]
    Set \( r_0 = \min\left\{1, \frac{\pi}{4\sqrt{K}} \right\}\)
    and define the regions 
    \[
    \Omega_k = \left\{(s, \exp X): |s| < r_0/2, \; |X|_{h_k(s)_\mathsf{e}} < r_0/2\right\}.
    \]
    As \(\Omega_k \subseteq 
    B_{\tilde{g}_k}(\tilde{p}_k, r_0)\), 
    \[
    \Vol_{\tilde{g}_k} B_{\tilde{g}_k}(\tilde{p}_k, r_0) \geq \Vol_{\tilde{g}_k} \Omega_k = \int_{-r_0/2}^{r_0/2}\Vol_{h_k(s)} \exp \left(B_{h_k(s)}(0, r_0/2)\right) ds
    = \omega_n\left( \frac{r_0}{2} \right)^n r_0. 
    \]
    From the curvature bound, the Cheeger--Gromov--Taylor injectivity radius estimate \cite[Thm. 4.7]{cheegergromovtaylor1982finite} gives 
    \[
    \inj(\tilde{p}_k, \tilde{g}_k) \geq \frac{r_0}{2} 
    \frac{1}{1+({V^{-K}_{2r_0}}/ \Vol_{\tilde{g}_k} B_{\tilde{g}_k}(\tilde{p}_k, r_0) )}
    \]
    where \(V^{-K}_{r}\) is the volume of a ball of radius \(r\) in the simply connected space form of constant curvature \(-K\). 
    This provides a lower bound 
    \[
    \inj(\tilde{p}_k, \tilde{g}_k) \geq
    \frac{r_0}{2}  
    \frac{1}{1+ ({V^{-K}_{2r_0}}/\omega_n\left(\frac{r_0}{2}\right)^n r_0 )} =: i(n, K)>0.
    \]
    By Hamilton's compactness theorem \cite{hamilton1996compactness} and Shi's derivative estimates \cite{shi1989deforming}, there is a subsequence such that 
    \begin{equation}\label{equ:pf-compact-universal-cover-converge}
        (\tilde{M}_k, \tilde{g}_k(t), \tilde{p}_k) \to (\tilde{M}_\infty, \tilde{g}_\infty(t), \tilde{p}_\infty)
    \end{equation}
    in the pointed smooth Cheeger--Gromov topology on compact time intervals. 
    Here \(\tilde{g}_\infty(t)\) is a complete normalised Ricci flow.
    \newline
    \indent
    Fix \(t = 0 \) and take \(\tilde{g}_k := \tilde{g}_k(0).\)
    On \(\tilde{M}_k\) there is a global right \(\mathsf{N}_k\)-action, and hence a Lie algebra of globally defined Killing fields of \((\tilde{M}_k, \tilde{g}_k)\). 
    For every \(k\), choose a vertical basis of Killing fields 
    \begin{equation}\label{equ:nk-Killing-basis}
         \mathfrak{n}_k := {\rm span}\{V_{1,k}^\#, \dots V_{n,k}^\#\}
    \end{equation}
    that are \(\tilde{g}_k\)-orthonormal at \(\tilde{p}_k\).
    The \(1\)-jet norm of a vector field is
    \[
    |U|^{*} := \left(|U|^2 + |\nabla U|^2 \right)^{1/2}.
    \]
    If \(U^\#\) is a vertical Killing field that is non-zero at \(p \in M\), then its \(1\)-jet norm satisfies the inequality
    \begin{equation}\label{equ:pf-compact-1-jet-bound}
        |U^\#|^{*}_p  \leq \left(1+ 2|L|^2_p +|A|^2_p + \frac{9}{4}|\mu|^2_p\right)^{1/2} |U^\#(p)|, 
    \end{equation}
    where \(L\) is the shape operator of the fibres, \(A\) is O'Neill's integrability tensor, and \(\mu\) is the bracket of the Lie algebra of Killing fields. 
    As \(A = 0\) for all \((\tilde{M}_k, \tilde{g}_k(t))\), to show that the \(1\)-jet norms \(|V_{i, k}^\#|^*_{\tilde{p}_k}\) are bounded, we must bound the shape operators and the structure constants 
    \[
    (\mu_k)_{ij}^m V_{m, k}^\# = [V_{i, k}^\#, V_{j, k}^\#].
    \]
    \indent
    Again, fix \(t = 0\) and take \(\tilde{g}_k : = \tilde{g}_k(0)\). 
    In the identification \eqref{equ:pf-compact-identify}, let \(\nabla^k\) be the Levi-Civita connections. 
    The shape operators of the fibres are
    \[
    L_k U = \nabla^k_U \partial_s
    \]
    where \(U\) is vertical. 
    The shape operator satisfies the Riccati equation 
    \begin{equation}\label{equ:pf-compact-Riccati}
        \nabla_{\partial_s}^k L_k + L_k^2 +  R_k(\cdot, \partial_s)\partial_s = 0. 
    \end{equation}
    Along a horizontal geodesic, take a parallel vertical unit vector field \(V\). 
    The function \(\varphi(s)= h_k(L_k V, V)\) is defined for all \(s \in \RR\) and satisfies 
    \[
    \frac{d \varphi}{ds} \leq K - \varphi^2. 
    \]
    If \(\varphi(s_0) < -\sqrt{K}\), comparison with \(\psi' = K - \psi^2\) forces blow-up in finite positive time, contradicting smoothness on all of \(\RR\).
    If \(\varphi(s_0) > \sqrt{K}\), then a contradiction follows by applying the previous argument to \(-\varphi(2s_0 - s)\).
    Hence, 
    \begin{equation}\label{equ:pf-compact-unit-vec-shape-bound}
    -\sqrt{K}\leq h_k(L_k V, V) \leq \sqrt{K}
    \end{equation}
    for all vertical unit vectors \(V\).
    The scalar curvatures of the nilpotent fibres are 
    \[
    \operatorname{scal}(h_k(s)) = -\frac{1}{4} |\mu_k|^2_{h_k(s)}.
    \]
    The traced Gauss equation and \eqref{equ:pf-compact-unit-vec-shape-bound} imply 
    \[
    |\operatorname{scal}(h_k(s)) | \leq n(n-1)K + |(\tr L_k)^2 - |L_k|^2| \leq 2n^2 K, 
    \]
    so the structure constants are bounded
    \begin{equation}\label{equ:pf-compact-bracket-bound}
        |\mu_k|^2  \leq 8n^2 K. 
    \end{equation}
    By \eqref{equ:pf-compact-1-jet-bound}, \eqref{equ:pf-compact-unit-vec-shape-bound}, and \eqref{equ:pf-compact-bracket-bound}, at \(\tilde{p}_k\) the \(1\)-jet norms of the Killing fields \(\{V_{i, k}^\#\}\) given in \eqref{equ:nk-Killing-basis}  are bounded by
    \begin{equation}\label{equ:pf-compact-1-jet-squeeze}
         1 \leq  \left|V_{i,k}^\#\right|^*_{\tilde{p}_k} \leq B(n, K) .
    \end{equation}
    Since the \(1\)-jet norms are bounded, 
    the sequences \((V_{i, k}^\#)_{k \in \NN}\) subconverge in the \(C^1\)-topology to \(\tilde{g}_\infty\)-Killing fields \(V_{i, \infty}^\#\) that are orthonormal at \(\tilde{p}_\infty\) \cite[p.~330]{heber1998noncompact}. 
    By the same argument as in \cite[Prop. 9.8]{bohmlaf2018},
    the brackets converge
    \[
    [V_{i,k}^\#, V_{j,k}^\#]
    \to [V_{i,\infty}^\#, V_{j,\infty}^\#]. 
    \]
    This means that
    \[
    \mathfrak{n}_\infty :=  {\rm span}\{V_{1,\infty}^\#, \dots, V_{n,\infty}^\#\}
    \]
    is an \(n\)-dimensional Lie algebra of Killing fields on \((\tilde{M}_\infty, \tilde{g}_\infty)\). 
    This Lie algebra is nilpotent since brackets of length \(n +1\) vanish for every \(\mathfrak{n}_k\), and this property passes to the limit.
    Since \((\tilde{M}_\infty, \tilde{g}_\infty)\) is complete, its Killing fields are complete. 
    Hence \(\mathfrak{n}_\infty\) integrates to a global isometric action on \((\tilde{M}_\infty, \tilde{g}_\infty)\). 
    \newline
    \indent
    Next, we pass the projections \(\tilde{M}_k\) onto \(\RR\) to the limit. 
    The shape operator satisfies 
    \[
     L_k = \frac{1}{2} h_k^{-1} \frac{d h_k}{ds}. 
    \]
    From the bound (\ref{equ:pf-compact-unit-vec-shape-bound})
    and by integrating, we find that the fibre metrics are uniformly equivalent on compact intervals
    \begin{equation}\label{equ:pf-compact-fibre-metric-uniform}
        e^{-2|s|\sqrt{K} }h_k(0) \leq h_k(s) \leq e^{2|s|\sqrt{K}} h_k(0). 
    \end{equation}  
    The Riccati equation \eqref{equ:pf-compact-Riccati} gives 
    \[
    |\nabla^k_{\partial_s} L_k| \leq (n+1)K. 
    \]
    Differentiating \eqref{equ:pf-compact-Riccati} and using Shi's estimates \cite{shi1989deforming} gives, inductively, 
    \begin{equation}\label{equ:pf-compact-derivative-shape-bound}
         |(\nabla^k_{\partial_s})^r L_k| \leq C_r. 
    \end{equation}
    Vertical derivatives are controlled by invariance, Koszul's formula, \eqref{equ:pf-compact-unit-vec-shape-bound}, and the bracket bound \eqref{equ:pf-compact-bracket-bound}. 
    Hence, all derivatives are uniformly bounded on compact \(s\)-intervals. 
    Let \(s_k: \tilde{M}_k \to \RR\) be the projections such that \(s_k(s, \mathsf{n}) = s\) via the identifications \eqref{equ:pf-compact-identify}.
    Its gradient and Hessian are 
    \[
    \nabla s_k = \partial_s, \qquad \operatorname{Hess}(s_k) = \tilde{g}_k(L_k \cdot, \cdot), 
    \]
    where \(L_k\) vanishes in the horizontal direction.
    Let \((\Phi_k)_{k\in\NN}\) be the exhausting pointed embeddings of the Cheeger--Gromov convergence \eqref{equ:pf-compact-universal-cover-converge}.     
    The estimates (\ref{equ:pf-compact-unit-vec-shape-bound}, \ref{equ:pf-compact-bracket-bound}, \ref{equ:pf-compact-derivative-shape-bound}) give bounds on compact sets for \(s_k \circ \Phi_k\) and all their derivatives. 
    Passing to a subsequence, \(s_k \circ \Phi_k\) converges to \(s_\infty\) smoothly on compact subsets, with 
    \[
    s_\infty(\tilde{p}_\infty) = 0, \qquad |\nabla s_\infty|_{\tilde{g}_\infty} = 1.
    \]
    The gradient field \(\nabla s_\infty\) is complete. 
    Let \(\theta_\infty^s\) be the flow of \(\nabla s_\infty\) and set \(\Sigma_\infty = s_\infty^{-1}(0)\). 
    The flow from \(\Sigma_\infty\)
    \begin{equation}\label{equ:pf-compact-hypersurface-flow}
        \RR \times \Sigma_\infty \to \tilde{M}_\infty, \qquad
    (s, x) \mapsto \theta^s_\infty(x)
    \end{equation}
    is a diffeomorphism. Indeed, \(\Sigma_\infty\) is connected and the inherited metric is complete. 
    Set \(\Sigma_k= s_k^{-1}(0)\) and let \(\theta^u_k\) be the local flow of \(\nabla^{\Phi_k^*\tilde{g}_k}(s_k \circ \Phi_k)\). 
    This gradient field is unit:
    \[
    \frac{d}{du} (s_k \circ \Phi_k)(\theta^u_k(x)) = |\nabla(s_k \circ \Phi_k)|^2 = 1. 
    \]
    By integrating, we get 
    \[
    (s_k \circ \Phi_k)(\theta^u_k(x)) = (s_k \circ \Phi_k)(x) + u. 
    \]
    Indeed, the map
    \[
    \hat{\Phi}_k(x) := \Phi_k \circ \theta_k^{-(s_k \circ \Phi_k)(x)}(x),
    \]
    defined on exhausting subsets of \(\Sigma_\infty\),
    maps \(\Sigma_\infty\) to \(\Sigma_k\). 
    The map \(x \mapsto  \theta_k^{-(s_k \circ \Phi_k)(x)}(x)\) converges smoothly on compact sets of \(\Sigma_\infty\) to the identity, so 
    \[
    \hat{\Phi}_k^*(\tilde{g}_k|_{\Sigma_k}) \to \tilde{g}_\infty|_{\Sigma_\infty}, 
    \qquad \hat{\Phi}_k(\tilde{p}_\infty) = \tilde{p}_k.
    \]
    These maps give smooth pointed Cheeger--Gromov convergence of the fibres. 
    Each \((\Sigma_k, \tilde{g}_k|_{\Sigma_k})\) is homogeneous and diffeomorphic to \(\RR^n\), so by \cite[Thm.~8.3]{bohmlaf2018} \(\Sigma_\infty\) is diffeomorphic to \(\RR^n\) too. 
    \newline
    \indent
    Recall that the Lie algebra of Killing fields \(\mathfrak{n}_\infty\) integrates to a global isometric right action by a connected and simply connected Lie group \(\tilde{\mathsf{N}}_\infty\). 
    As \(V_{i, k}^\#(s_k) = 0\), the limit satisfies \(V_{i, \infty}^\#(s_\infty) = 0\), so the \(\tilde{\mathsf{N}}_\infty\)-action preserves the level sets of \(s_\infty\).
    Let \(\mathsf{N}_\infty\) be the quotient of \(\tilde{\mathsf{N}}_\infty\) by the effective kernel. 
    As
    \[
    T_{\tilde{p}_\infty} \Sigma_\infty = 
    \operatorname{span}\{V_{1,\infty}^\#(\tilde{p}_\infty), \dots, V_{n,\infty}^\#(\tilde{p}_\infty)\},
    \]
    the effective kernel is discrete. 
    Hence, \({\rm Lie}(\mathsf{N}_\infty) = \mathfrak{n}_\infty\), so \(\mathsf{N}_\infty\) is nilpotent and \(\dim \mathsf{N}_\infty = n\).
    \newline
    \indent
    The orbit map
    \begin{equation}\label{equ:pf-compact-orbit-map}
        \mathsf{N}_\infty \to \Sigma_\infty, \qquad \mathsf{n} \mapsto \tilde{p}_\infty \cdot \mathsf{n}
    \end{equation}
    is a local diffeomorphism because its differential at the identity is the fundamental vector field point map \(\mathfrak{n}_\infty\to T_{\tilde{p}_\infty}\Sigma_\infty\) such that \(X \mapsto X^\#_{\tilde{p}_\infty}\). 
    Pull back the metric \(\tilde{g}_\infty |_{\Sigma_\infty}\) by this map to get a right-invariant metric \(h_\infty(0)\) on \(\mathsf{N}_\infty\). 
    Indeed, the orbit map is a local isometry from a complete connected Riemannian manifold, and thus is a covering. 
    This means that the \(\mathsf{N}_\infty\)-action on \(\Sigma_\infty\) is transitive. 
    The action is also free. 
    This is because stabilisers act orthogonally on the Lie algebra. Since \(\mathsf{N}_\infty\) is connected and nilpotent, the Lie exponential map is surjective. 
    If \(\mathsf{n} \in \mathsf{N}_\infty\) is a stabiliser, then there exists an \(X \in \mathfrak{n}_\infty\) such that \(\mathsf{n} = \exp X\). 
    Hence, \(\operatorname{Ad}(\mathsf{n}) = \exp(\operatorname{ad}X)\) is orthogonal and unipotent, so it is the identity. 
    This means that \(\mathsf{n}\) is central. 
    By transitivity, \(\mathsf{n}\) fixes all of \(\Sigma_\infty\). 
    As \([V_{i,\infty}^\#, \nabla s_\infty] = 0\), the \(\mathsf{N}_\infty\)-actions commute with \(\theta_\infty^s\), so \(\mathsf{n}\) fixes all of \(\tilde{M}_\infty\). 
    But the \(\mathsf{N}_\infty\) action is effective, so \(\mathsf{n} = \mathsf{e}\).
    \newline
    \indent
    As the action on \(\Sigma_\infty\) is free and transitive, and \(\Sigma_\infty\) is diffeomorphic to \(\RR^n\), we get \(\mathsf{N}_\infty\) is simply-connected.
    From \eqref{equ:pf-compact-hypersurface-flow} and \eqref{equ:pf-compact-orbit-map}, we get the diffeomorphism 
    \[
    \RR \times \mathsf{N}_\infty \to \tilde{M}_\infty, \qquad (s, \mathsf{n}) \mapsto \theta^s_\infty(\tilde{p}_\infty\cdot\mathsf{n}), 
    \]
    under which
    \[
    (\tilde{M}_\infty, \tilde{g}_\infty, \tilde{p}_\infty) \cong (\RR \times \mathsf{N}_\infty, ds^2 +h_\infty(s), (0, \mathsf{e})).
    \]
    Here \(h_\infty(s)\) is right-invariant, and the action is free and transitive on the fibres. 
    The limit is a principal \(\mathsf{N}_\infty\)-bundle over \(\RR\) equipped with an invariant metric.
    This proves the first assertion. 
    \newline
    \indent
    The vertical Killing fields of \eqref{equ:nk-Killing-basis} are independent of time and Killing along the flow \(\tilde{g}_k(t)\). 
    The smooth convergence on compact sets gives the same action along the limit flow \(\tilde{g}_\infty(t)\). 
    \newline
    \newline
    \indent
    \textbf{Part (b).}
    Next, assume that the base diameter is bounded below by some \(c >0\). 
    At time \(t = 0\), set 
    \[
    \ell_k = 2\diam(S^1, \check{g}_k(0)), \qquad 2c\leq \ell_k \leq 2d.
    \]
    After passing to a subsequence, \(\ell_k\) converges to \(\ell_\infty \in [2c, 2d]\). 
    Take \(\tau_k\) to be the deck generators of \(\tilde{M}_k \to M_k\) such that in \eqref{equ:pf-compact-identify} coordinates
    \[
    \tau_k(s, \mathsf{n}) = (s +\ell_k, \varphi_k^{-1}(\mathsf{n})), \qquad h_k(s+\ell_k) = \varphi_k^*h_k(s)
    \]
    where \(\varphi_k \in {\rm Aut}(\mathsf{N}_k)\ltimes(\mathsf{N}_k)_L\) is the holonomy. 
    Set \(\mathsf{h}_k = \varphi_k(\mathsf{e})\). Then the maps
    \[
    \sigma_k := \tau_k \circ R_{\mathsf{h}_k}
    \]
    are isometries of \(\tilde{g}_k(t)\) and at \(t = 0\)
    \[
    \sigma_k(0, \mathsf{e}) = (\ell_k, \mathsf{e}), 
    \qquad \operatorname{dist}_{\tilde{g}_k} ((0,\mathsf{e}), \sigma_k(0, \mathsf{e})) = \ell_k \leq 2d.
    \]
    Compactness of isometries with bounded basepoint displacement (cf. \cite[Lem. B.3]{topping2010pseudo}) therefore gives, on a further subsequence 
    \[
    \Phi_k^{-1} \circ \sigma_k \circ \Phi_k \to \sigma_\infty 
    \]
    smoothly on compact subsets. 
    As \(s_\infty \circ \sigma_\infty = s_\infty +\ell_\infty\) and \(\ell_\infty >0\), the action of \(\langle \sigma_\infty \rangle\) is free and properly discontinuous. 
    Also, \(\sigma_\infty\) commutes with the gradient flow of \(s_\infty\), so it is of the form 
    \[
    \sigma_\infty(s, \mathsf{n}) = (s +\ell_\infty, \varphi_\infty(\mathsf{n}))    
    \]
    where \(\varphi_\infty: (\mathsf{N}_\infty, h_\infty(0)) \to (\mathsf{N}_\infty, h_\infty(\ell_\infty))\) is an isometry. 
    Isometries between right-invariant metrics on connected nilpotent Lie groups are affine \cite[Thm. 1.2]{kivioja2017isometries}. 
    Then, \(\varphi_\infty \in {\rm Aut}(\mathsf{N}_\infty) \ltimes (\mathsf{N}_\infty)_L\) and 
    \[
    (\RR \times \mathsf{N}_\infty)/\langle \sigma_\infty \rangle \to \RR /\ell_\infty \ZZ 
    \]
    is a twisted principal \(\mathsf{N}_\infty\)-bundle over \(S^1\). 
    The metrics \(\tilde{g}_\infty(t)\) descend to invariant metrics solving the normalised Ricci flow. 
    \newline
    \newline
    \indent
    \textbf{Part (c).}
    It remains to show convergence of \((M_k, g_k(t), p_k)\) under the assumptions 
    \[
     0 < c\leq \diam(S^1, \check{g}_k(t)),\qquad
     \operatorname{dist}_{\tilde{g}_k(0)}(\tilde{p}_k, \tau_k(\tilde{p}_k)) \leq l,
    \]
    for some \(c, l >0\). 
    First, we prove an injectivity radius estimate. 
    Set \(g_k(0) = g_k\) and
    \[
    r_1 = \min\left\{c, \frac{\pi}{4\sqrt{K}}\right\}.
    \]
    Take horizontal unit-speed geodesics \(\gamma_k:(-\frac{r_1}{2}, \frac{r_1}{2}) \to M_k\) with \(\gamma_k(0) = p_k\). 
    Define the regions
    \[
    \Omega_k = \{\gamma_k(s) \cdot \exp X:|s|< r_1/2, \; |X|_{h_k(s)_\mathsf{e}}< r_1/2\}. 
    \]
    As before, 
    \(\Omega_k \subseteq 
    B_{{g}_k}({p}_k, r_1)\), and so
    \[
    \Vol_{{g}_k} B_{{g}_k}({p}_k, r_1) \geq \Vol_{{g}_k} \Omega_k = \int_{-r_1/2}^{r_1/2}\Vol_{h_k(s)} \exp \left(B_{h_k(s)}(0, r_1/2)\right) ds
    = \omega_n\left( \frac{r_1}{2} \right)^n r_1. 
    \]
    The Cheeger--Gromov--Taylor injectivity radius estimate \cite[Thm. 4.7]{cheegergromovtaylor1982finite} and the curvature bound yield
    \[
    \inj({p}_k, {g}_k) \geq
    \frac{r_1}{2}  
    \frac{1}{1+ ({V^{-K}_{2r_1}}/\omega_n\left(\frac{r_1}{2}\right)^n r_1 )} =: i(n,c, K)>0.
    \]
    By Hamilton's compactness Theorem \cite{hamilton1996compactness}, after passing to a subsequence
    \begin{equation}\label{equ:pf-compact-tpb-converge}
        (M_k, g_k(t), p_k) \to (M_\infty, g_\infty(t), p_\infty)
    \end{equation}
    in the pointed smooth Cheeger--Gromov topology. 
    Let \((\Psi_k)_{k\in\NN}\) be the exhausting pointed embeddings of this Cheeger--Gromov convergence. 
    \newline
    \indent
    As we have assumed that the displacement of \(\tilde{p}_k\) by the deck generators \(\tau_k\) is bounded, compactness of isometries (cf. \cite[Lem. B.3]{topping2010pseudo}) gives 
    \[
    \Phi_k^{-1} \circ \tau_k \circ \Phi_k \to \tau_\infty
    \]
    smoothly on compact subsets. 
    By the same argument as for \(\sigma_\infty\), the group \(\langle\tau_\infty\rangle\) acts freely and properly discontinuously, and 
    \[
    \tau_\infty(s, \mathsf{n}) = (s+\ell_\infty, \phi_\infty(\mathsf{n}))
    \]
    for some affine \(\phi_\infty \in {\rm Aut}(\mathsf{N}_\infty) \ltimes (\mathsf{N}_\infty)_L\). 
    To identify the limit in \eqref{equ:pf-compact-tpb-converge} with this quotient, let \(q_k : \tilde{M}_k \to M_k\) be the covering maps. 
    Then the maps 
    \[
    \Psi_k^{-1} \circ q_k \circ \Phi_k
    \]
    are defined on exhausting subsets of \(\tilde{M}_\infty\) and send \(\tilde{p}_\infty\) to \(p_\infty\). 
    Their images of compact subsets remain in bounded metric balls since \(q_k\) does not increase distances.
    As these maps are local isometries for the pulled-back metrics, their derivatives are controlled by the smoothly converging metrics. 
    Passing to a further subsequence gives a smooth limit 
    \[
    q_\infty : \tilde{M}_\infty \to M_\infty, \qquad q_\infty^* g_\infty(t) = \tilde{g}_\infty(t)
    \]
    for compact time intervals. 
    Indeed, \(q_\infty\) is a covering and \(q_\infty \circ \tau_\infty = q_\infty\). 
    Finally, we show that the fibres of \(q_\infty\) are the orbits of \(\tau_\infty\). 
    Suppose that \(q_\infty(x) = q_\infty(y)\). 
    Choose sufficiently small neighbourhoods of \(x\) and \(y\) such that \(q_\infty\) is a diffeomorphism. 
    By smooth convergence and the inverse function theorem, there exist sequences of points \(x_k, y_k \in \tilde{M}_k\) such that 
    \[
    \Phi_k^{-1}(x_k) \to x, \qquad \Phi_k^{-1}(y_k) \to y, \qquad q_k(x_k) = q_k(y_k). 
    \]
    Indeed, \(y_k = \tau^{m_k}_k(x_k)\) for some \(m_k \in \ZZ\). 
    These points lie in \(B_{\tilde{g}_k}(\tilde{p}_k, R)\) for fixed \(R\). 
    \newline
    \indent
    Since the fibres are simply connected, the bundle projections \(M_k \to S^1\) induce an isomorphism of fundamental groups. 
    Every curve joining \(x_k\) to \(\tau^{m_k}_k(x_k)\) projects to a loop in \(S^1\) with winding number \(m_k\). 
    Indeed, 
    \[
    |m_k|\ell_k \leq \operatorname{dist}_{\tilde{g}_k}(x_k, y_k)\leq 2R, 
    \]
    and since \(\ell_k \geq 2c\), we have \(|m_k| \leq R/c\). 
    After passing to a subsequence, \(m_k = m\) and \(y = \tau_\infty^m(x)\), so 
    \[
    (\tilde{M}_\infty, \tilde{g}_\infty(t))/\langle \tau_\infty\rangle \cong (M_\infty, g_\infty(t)).
    \]
\end{proof}

\begin{prop}\label{prop:adjoint-holonomy-compactness}
    Fix \(n \in \mathbb{N}\) and \(K,c,  d >0\). 
    Let \(I\subseteq \RR\) be an open interval containing \(0\). 
    Let \(\pi_k:M_k \to S^1\) be a sequence of twisted principal \(\mathsf{N}_k\)-bundles where \(\mathsf{N}_k\) is nilpotent, simply-connected and \(\dim \mathsf{N}_k = n\). 
    For every \(k \in \NN\), let \(g_k(t)\) be a normalised Ricci flow of invariant metrics on \(M_k\) defined for \(t \in I\) such that 
    \[
    |\Rm_{g_k(t)}|_{g_k(t)} \leq K, \qquad 0< c\leq \diam(S^1, \check{g}_k(t)) \leq d. 
    \]
    Take a sequence of points \(p_k \in M_k\) and their lifts to the universal covers \(\tilde{p}_k \in \tilde{M}_k\).
    At \(t = 0\), let \(P_k \in {\rm Aut}(\mathfrak{n}_k)\) be the adjoint holonomy of \(\nabla^{{\mathrm{ad}}}_k\) based at \(\pi_k(p_k)\). 
    Let \((\tilde{M}_\infty, \tilde{g}_\infty(t), \tilde{p}_\infty)\) be a subsequential limit that descends to \(\pi_\infty:M_\infty \to S^1\) as in \autoref{thm:inv-RF-compactness}. 
    Then, there exist linear isometries between the fibres of the adjoint bundle \(\rho_k : \mathfrak{n}_\infty \to \mathfrak{n}_k\) such that 
    \[
    \rho_k^{-1} P_k \rho_k \to P_\infty
    \]
    where \(P_\infty\) is the holonomy of \(\nabla^{\rm ad}_\infty\) based at \(\pi_\infty(p_\infty)\) at time \(t = 0\).     
\end{prop}

\begin{proof}
    Let \(\{V_{i, k}^\#\}\) be the converging orthonormal bases \eqref{equ:nk-Killing-basis} in the proof of \autoref{thm:inv-RF-compactness}.
    Define \(\rho_k\) by \(\rho_k V_{i, \infty}^\# = V^\#_{i, k}\).   
    In the horizontal trivialisation based at \(\tilde{p}_k\), write 
    \[
    \tau_k(s, \mathsf{n}) =(s + \ell_k, \varphi_k^{-1}(\mathsf{n})) 
    \]
    where \(\ell_k\) is the length of \((S^1, \check{g}_k(0))\) and \(\varphi_k\) is the twisted principal bundle holonomy. 
    Define \(\psi_k = R_{\varphi_k(\mathsf{e})}^{-1} \circ \varphi_k\) and \(\sigma_k = \tau_k \circ R_{\varphi_k(\mathsf{e})}\). 
    The adjoint holonomy is \(P_k =  d(\psi_k)_\mathsf{e}\) and 
    \[
    \sigma_k(s, \mathsf{n}) = (s+\ell_k, \psi_k^{-1}(\mathsf{n})). 
    \]
    As in the proof of \autoref{thm:inv-RF-compactness}, after passing to a subsequence, \(\ell_k \to \ell_\infty >0\) and, in the Cheeger--Gromov maps, \( \sigma_k \to \sigma_\infty\) where 
    \[
    \sigma_\infty(s, \mathsf{n}) = (s +\ell_\infty, \psi_\infty^{-1}(\mathsf{n}))
    \]
    for some \(\psi_\infty \in {\rm Aut}(\mathsf{N}_\infty)\). 
    In the Killing basis
    \[
    (\sigma_k)_* V_{i, k}^\# = \sum_{j = 1}^n  (P^{-1}_k)_i^j V_{j, k}^\#. 
    \]
    By the convergence of \(\sigma_k\) and the vertical Killing fields, we have 
    \[
    \rho_k^{-1} P_k^{-1} \rho_k \to d(\psi_\infty^{-1})_\mathsf{e}
    \]
    The limit flow descends to \(M_\infty = (\RR \times \mathsf{N}_\infty)/\langle \sigma_\infty \rangle\) that has bundle holonomy \(\psi_\infty\) and adjoint holonomy \(P_\infty = d(\psi_\infty)_\mathsf{e}\). 
\end{proof}

\subsection{Subconvergence to solvsolitons}
Recall that a solvsoliton is a left-invariant Ricci soliton metric on a solvable Lie group \(\mathsf{S}\) such that when evaluated at the identity, the inner product on the Lie algebra satisfies
\begin{equation}\label{equ:alg-soliton}
    \Ric_g = c\operatorname{Id} +D,
\end{equation}
for some \(c \in \RR\) and some derivation of the Lie algebra \(D\).  
Similarly, a nilsoliton is a left-invariant Ricci soliton metric on a nilpotent Lie group \(\mathsf{N}\) that satisfies \eqref{equ:alg-soliton}. 
We call right-invariant metrics solvsolitons or nilsolitons too if they satisfy \eqref{equ:alg-soliton}, since inversion identifies these metrics with left-invariant metrics with the same inner product at the identity. 
\newline
\indent
Let \(\mathsf{N}\) be a simply-connected nilpotent Lie group.
On twisted principal  \(\mathsf{N}\)-bundles \(M \to S^1\), consider the normalised Ricci flow of invariant metrics in the following gauge
 \begin{equation}\label{equ:normalised-Ricci-flow-gauged}
     \frac{\partial g}{\partial t} = -2\ric_g  - g + \mathcal{L}_N g, \qquad g(0) = g_0, 
 \end{equation}
where \(N = -(\tr L)\bar X\) is the mean curvature vector of the fibres. 
Given a sequence \(t_k \to \infty\), fix a point \(p \in M\) and take the sequence of time-translated flows \(g_k(t) := g(t_k +t)\). 
By \autoref{thm:inv-RF-compactness}, the flows lifted to the universal cover subconverge to an invariant solution to \eqref{equ:normalised-Ricci-flow-gauged} on a principal \(\mathsf{N}_\infty\)-bundle \(\tilde{M}_\infty \to \RR\) for a simply-connected nilpotent Lie group \(\mathsf{N}_\infty\). 
Under the additional assumption that the diameter of \((S^1, \check{g}_k(t))\) is bounded below, we extract properties of the limit \((\tilde{M}_\infty, \tilde{g}_\infty(t))\) using the
Gindi--Streets functional \cite[Def. 4.1]{gindi2023four}
\[
\mathcal{I}(g) = \int_{S^1} |L|^2 + \frac{1}{2} ds.
\]
This functional is invariant under twisted principal bundle isomorphisms.  

\begin{lemma}\label{lemma:Gindi-Streets-monotone}
    Let \(M \to S^1\) be a twisted principal \(\mathsf{N}\)-bundle where \(\mathsf{N}\) is a nilpotent Lie group. 
    If \(g(t)\) is an invariant normalised Ricci flow on \(M\), then
    \[
    \frac{d}{dt} \mathcal{I}(g(t))=-\int_{S^1}|\delta_{[\cdot, \cdot]}(L)|^2 + \left(|L|^2  - \frac{1}{2}\right)^2 + 4|\hat{\nabla}^*L |^2 + 2|\nabla L|^2 \;ds, 
    \]
    where \([\cdot, \cdot]\) is the Lie bracket on the adjoint bundle.     
\end{lemma}

\begin{proof}
    By \autoref{prop:Ricci-curv-dim1}, under the normalised flow in the mean curvature gauge \eqref{equ:normalised-Ricci-flow-gauged} the volume form of the base evolves by 
    \[
    \partial_t ds  = \left(|L|^2  - \frac{1}{2}\right) ds. 
    \] 
    By equation \eqref{equ:heat-shape-operator} in the proof of \autoref{prop:shape-operator-estimate}, 
    \begin{align*}
        \left(\pdv{}{t} - \Delta_{S^1}\right) |L|^2
        &= -2h(\nabla \hat\Ric_g, L)
        -2|L|^4 + |L|^2  -4h(\hat{\nabla}^*L - L H, \hat{\nabla}^* L) - 2|\nabla L|^2
        \\
        &= -|\delta_{[\cdot, \cdot]}(L)|^2  -2|L|^2\left(|L|^2  - \frac{1}{2}\right)  - 4|\hat{\nabla}^*L |^2 - 2|\nabla L|^2, 
    \end{align*}
    where the second line follows from nilpotence and \autoref{prop:1st-2nd-variation-moment-map}.
    Combining these, we get
    \[
    \frac{d}{dt} \mathcal{I}(g(t)) = -\int_{S^1}|\delta_{[\cdot, \cdot]}(L)|^2 + \left(|L|^2  - \frac{1}{2}\right)^2 + 4|\hat{\nabla}^*L |^2 + 2|\nabla L|^2 ds. 
    \]
\end{proof}

\begin{lemma}\label{lemma:four-properties}
    Let \(M \to S^1\) be a twisted principal \(\mathsf{N}\)-bundle where \(\mathsf{N}\) is a nilpotent Lie group. 
    Let \(g(t)\) be an invariant normalised Ricci flow on \(M\). 
    Then, the length of the base \(\ell(t) := 2\diam(S^1, \check{g}(t))\) converges to \(\ell_\infty \geq 0\) as \(t \to \infty\) where 
    \[
     \ell_\infty =  \lim_{t \to \infty} \mathcal{I}(g(t)). 
    \]
    Furthermore, if \(\ell_\infty> 0\), then every subsequential limit of time-translated flows obtained in \autoref{thm:inv-RF-compactness} satisfies 
    \[
    \delta_{[\cdot, \cdot]}( L_\infty )=0, 
    \qquad |L_\infty|^2 = \frac{1}{2}, 
    \qquad {\hat\nabla}^*L_\infty = 0, \qquad \nabla L_\infty = 0. 
    \]
    along the limit flow. 
\end{lemma}

\begin{proof}
    By \autoref{lemma:Gindi-Streets-monotone}, 
    the Gindi--Streets functional \(\mathcal{I}(g(t))\) is monotone and nonnegative, and therefore converges to some \(\mathcal{I}_\infty \geq 0\). 
    Since 
    \[
    \frac{d \ell}{dt} = \int_{S^1} |L|^2 - \frac{1}{2}\, ds = \mathcal{I}(g(t)) - \ell(t),
    \]
    we get 
    \[
    \ell(t) = e^{-t} \ell(0) + \int_0^t e^{s-t} \mathcal{I}(g(s))\, ds.  
    \]
    Hence \(\ell(t) \to \mathcal{I}_\infty\) as \(t \to \infty\).
    Now assume that \(\ell_\infty >0\). 
    Take a sequence of time-translated normalised Ricci flows of invariant metrics \(g(t + t_k)\) that converge as in \autoref{thm:inv-RF-compactness}. 
    For every fixed \(a< b\), 
    \begin{equation}\label{equ:functional-integral-derivative}
         \mathcal{I}(g(t_k + a)) - \mathcal{I}(g(t_k +b))= \int_{t_k+a}^{t_k+b}\int_{S^1}|\delta_{[\cdot, \cdot]}(L)|^2 + \left(|L|^2  - \frac{1}{2}\right)^2 + 4|\hat{\nabla}^*L |^2 + 2|\nabla L|^2 \;ds\,dt 
    \end{equation}
    converges to zero as \(k \to \infty\).
    Assume that one of the terms 
    \[
    \delta_{[\cdot, \cdot]}( L_\infty ), 
    \qquad |L_\infty|^2 - \frac{1}{2}, 
    \qquad {\hat\nabla}^*L_\infty , \qquad \nabla L_\infty
    \]
    is non-zero at some point and time of the limit flow. 
    Then there is a neighbourhood and a time interval for which the corresponding term in the integral \eqref{equ:functional-integral-derivative} is bounded below by a positive constant for large \(k\). 
    This is a contradiction. 
\end{proof}

After we take a subsequential limit \((\tilde{M}_\infty, \tilde{g}_\infty(t))\), we run the normalised bracket flow. 
For the normalised Ricci flow of invariant metrics with gauge \eqref{equ:normalised-Ricci-flow-gauged}, the corresponding bracket flow solves the following system: 
\begin{align}
     \frac{\partial[\cdot, \cdot]_\mu}{\partial t} &= \delta_{[\cdot, \cdot]_\mu} \left( \hat{\Ric}_{\mu} -(\nabla L)+\frac{1}{2}\operatorname{Id}\right),
     \label{equ:bf-normal-dim-1-mu}
     \\
     \frac{\partial\omega_\mu}{\partial t}(\bar{X}) &= 2\omega_\mu ({\hat\nabla}^*L),
     \label{equ:bf-normal-dim-1-omega}
     \\
     \frac{\partial\check{g}_\mu}{\partial t} &=  (2|L|^2 -1) \check{g}_\mu. 
     \label{equ:bf-normal-dim-1-g^}
\end{align}
Just as in the proof of \autoref{thm:B}, one obtains an equivalence of the normalised Ricci flow and normalised bracket flow by solving the fibrewise ODEs 
\begin{enumerate}[(i)]
    \item \(
    \frac{\partial }{\partial t} \Phi = - \Phi \iota \left( \hat\Ric_g - \nabla L + \frac{1}{2}\operatorname{Id} \right) \omega, \qquad \Phi(0) = \operatorname{Id},\)
    \item 
    \(\frac{\partial }{\partial t} \Phi = -\iota \left( \hat\Ric_{\mu} -\nabla L + \frac{1}{2}\operatorname{Id}\right)\omega_\mu \Phi, \qquad \Phi(0) = \operatorname{Id}.\)
\end{enumerate}
Along the homogeneous normalised bracket flow, if \(\mu(t)\) is nilpotent, then the following quantity is monotone \cite[Prop 3.11]{lauret2013RicciFlowOf}  
\[
R_\mu + |\Ric_\mu|^2. 
\]
Inspired by this, we integrate this over the base to get the following functional
    \[
    \mathcal{F}(\mu, h, \check{g}) = \int_{S^1} \hat{R}_\mu + |\hat{\Ric}_\mu|^2\, ds.
    \]
This functional is invariant under principal bundle isomorphisms. 

\begin{lemma}\label{lemma:Lauret-monotone}
    Let \(M \to S^1\) be a twisted principal \(\mathsf{N}\)-bundle where \(\mathsf{N}\) is a nilpotent Lie group. 
    Fix a vector bundle metric \(h\) on \(\adM\) and then let \(([\cdot, \cdot]_\mu, \omega_\mu, \check{g}_\mu)(t)\) be a normalised bracket flow on \(M\) solving (\ref{equ:bf-normal-dim-1-mu}--\ref{equ:bf-normal-dim-1-g^}).  
    Suppose that the initial condition has the following properties 
    \[
    \delta_{\mu}(L)=0, 
    \qquad |L|^2 = \frac{1}{2}, 
    \qquad\hat{\nabla}^*L = 0, \qquad \nabla L = 0. 
    \]
    Along this flow, the functional 
    \[
    \mathcal{F}(\mu(t), h, \check{g}(t)) = \int_{S^1} \hat{R}_\mu + |\hat{\Ric}_\mu|^2\, ds
    \]
    is monotonically decreasing, and 
    \[
    \lim_{t \to \infty}
    \int_{S^1} \left|  \delta_{\mu}\left(\hat{\Ric}_{\mu} + \frac{1}{2} \operatorname{Id}\right)\right|^2ds = 0. 
    \]
\end{lemma}

\begin{proof}
    Given the initial condition, the normalised bracket flow equations in the mean curvature gauge (\ref{equ:bf-normal-dim-1-mu}--\ref{equ:bf-normal-dim-1-g^}) are  
    \[
    \frac{\partial[\cdot, \cdot]_\mu}{\partial t} = \delta_{[\cdot, \cdot]_\mu} \left( \hat{\Ric}_{\mu} +\frac{1}{2}\operatorname{Id}\right),
    \qquad
    \frac{\partial\omega_\mu}{\partial t} = 0, \qquad \frac{\partial\check{g}_\mu}{\partial t} = 0.
    \]
    As the metric \(\check{g}\) on the base and the metric on the adjoint bundle \(h\) are both constant, we write 
    \(\mathcal{F}(\mu) = \mathcal{F}(\mu, h_0, \check{g})\).
    Since \(\mathsf{N}\) is nilpotent, the Ricci curvature of the fibres is \(\hat{\Ric}_\mu = M_\mu\) and scalar curvature is \(\hat{R}_\mu = -\frac{1}{4}|\mu|^2\).
    By \cite[Lem. 6.2]{lauret2011nilmanifolds},
    \[
    \partial_t \hat{\Ric}_\mu = -\frac{1}{2}\Delta_\mu\left(\hat{\Ric}_\mu + \frac{1}{2} \operatorname{Id}\right), 
    \]
    where \(\Delta_\mu = S\circ \delta_\mu^*\delta_\mu\), and \(S(\alpha) = \frac{1}{2}(\alpha + \alpha^*)\), and \(\delta_\mu^*\) is the adjoint of \(\delta_\mu\). 
    The time derivative of the functional is
    \[
    \frac{d}{dt} \mathcal{F}(\mu(t)) 
    =
    \int_{S^1} \left\langle\partial_t \hat\Ric_\mu, \operatorname{Id} + 2\hat\Ric_\mu \right\rangle
    ds 
    =-  \int_{S^1} \left|  \delta_{\mu}\left(\hat{\Ric}_{\mu} + \frac{1}{2} \operatorname{Id}\right)\right|^2ds.
    \]
    Suppose for contradiction that there exists a sequence \(t_k \to \infty\) and an \(\varepsilon >0\) such that 
     \[
    -  \int_{S^1} \left|  \delta_{\mu(t_k)}\left(\hat{\Ric}_{\mu(t_k)} + \frac{1}{2} \operatorname{Id}\right)\right|^2ds= \frac{d }{d t}\Big|_{t = t_k} \mathcal{F}(\mu(t))   
    \leq -\varepsilon <0. 
    \]
    The second time derivative of the bracket on \(\adM\) is 
    \[
    \frac{\partial^2\mu}{\partial t^2} =  \frac{1}{2} \frac{\partial \mu}{\partial t} + p(\mu),
    \]
    where \(p(\mu)\) is some polynomial of \(\mu\). 
    By \autoref{prop:bracket-bound}, under the normalised bracket flow \(|\mu| \leq C\) for some constant \(C>0\).   
    Then, for some constant \(C>0\), 
    \[
    \Big| \frac{\partial^2\mu}{\partial t^2}\Big| \leq \frac{1}{2}\Big| \frac{\partial \mu}{\partial t}\Big| +C,
    \]
    and so 
    \[
    \frac{d^2}{d t^2} \mathcal{F}(\mu(t)) =-2 \int_{S^1} \left\langle  \frac{\partial^2\mu}{\partial t^2},  \frac{\partial\mu}{\partial t} \right\rangle ds \leq   \frac{d}{d t} \mathcal{F}(\mu(t))
    +C\ell,
    \]
    where \(\ell\) is the length of \((S^1, \check{g})\).
    By the maximum principle, 
    \[
    \frac{d}{d t} \mathcal{F}(\mu(t))\leq 
    e^{t-t_k}\left(C\ell+\frac{d}{d t}\Big|_{t= t_k} \mathcal{F}(\mu(t)) \right) - C\ell 
    \leq e^{t-t_k}\left(C\ell- \varepsilon\right) - C\ell
    \]
    for \(t\geq t_k\). 
    There exists a \(\delta >0\) such that  
    \[
    \frac{d}{d t} \mathcal{F}(\mu(t)) < -\frac{\varepsilon}{2}
    \]
    for \(t- t_k< \delta\).
    Pass to a subsequence with \(t_{k+1} >t_k +\delta\) for all \(k\).  
    By Cauchy--Schwarz 
    \[
    \frac{1}{n}\left(\hat{R}_\mu+ \frac{n}{2}\right)^2 - \frac{n}{4}
    = \hat{R}_\mu + \frac{1}{n}\hat{R}_\mu^2 
    \leq
    \hat{R}_\mu + |\hat{\Ric}_\mu|^2,
    \]
    so \(-\frac{n}{4}\ell\leq \mathcal{F}(\mu(t))\).
    This provides the following contradiction
    \[
    \sum_{k = 1}^\infty \frac{\delta \varepsilon}{2}
    \leq -\sum_{k = 1}^\infty \int_{t_k}^{t_k +\delta}\frac{d}{d t} \mathcal{F}(\mu(t))\, dt
    \leq  - \int_0^\infty \frac{d}{dt}\mathcal{F}(\mu(t))\, dt
    \leq \mathcal{F}(\mu(0)) + \frac{n}{4} \ell< \infty. 
    \]    
\end{proof}

\begin{proof}[Proof of \autoref{thm:C}]
    By \autoref{thm:inv-RF-compactness}, the time-translated normalised Ricci flows \(g_k(t) := g(t+ t _k)\) subconverge to an invariant normalised Ricci flow on a principal \(\mathsf{N}_\infty\)-bundle over \(\RR\),
    \[
    (\tilde{M}_\infty, \tilde{g}_\infty(t), \tilde{p}_\infty), 
    \]
    where \(\mathsf{N}_\infty\) is a simply connected nilpotent Lie group. 
    Let \(\ell(t) = 2\diam(S^1, \check{g}(t))\).
    By \autoref{lemma:four-properties}, \(\ell(t) \to \ell_\infty \geq 0\). 
    First, if \(\ell_\infty = 0\), then \((\tilde{M}_\infty, \tilde{g}_\infty(t))\) is homogeneous, and by \cite[Thm. A]{bohmlaf2018} and \cite[Coro. B]{bohmlafuente2023noncompact}, we get our desired conclusion. 
    \newline
    \indent
    Next, if \(\ell_\infty > 0\), then by \autoref{lemma:four-properties}, there exists a constant \(d>0\) such that 
    \[
    0< d \leq \diam(S^1, \check{g}(t))
    \]
    for all \(t \geq 0\). 
    By the same lemma, \((\tilde{M}_\infty, \tilde{g}_\infty(t))\) exhibits the following properties
    \[
    \delta_{[\cdot, \cdot]}( L)=0, 
    \qquad |L|^2 = \frac{1}{2}, 
    \qquad {\hat\nabla}^*L = 0, \qquad \nabla L = 0, 
    \]
    where \([\cdot, \cdot]\) is the Lie bracket on the adjoint bundle. 
    Let \(([\cdot, \cdot]_\mu, \omega_\mu, \check{g}_\mu)(t)\) be the bracket flow corresponding to the normalised Ricci flow \(\tilde{g}_\infty(t)\). 
    By the above properties, the normalised bracket flow equations in the mean curvature gauge (\ref{equ:bf-normal-dim-1-mu}--\ref{equ:bf-normal-dim-1-g^}) are  
    \[
    \frac{\partial[\cdot, \cdot]_\mu}{\partial t} = \delta_{[\cdot, \cdot]_\mu} \Big( \hat{\Ric}_{\mu} +\frac{1}{2}\operatorname{Id}\Big),
    \qquad
    \frac{\partial\omega_\mu}{\partial t} = 0, \qquad \frac{\partial\check{g}_\mu}{\partial t} = 0.
    \]
    As \(\RR\) is contractible, \(\tilde{M}_\infty\) is a trivial principal bundle. 
    By \autoref{prop:local-principal-bundle-realise}, the bracket flow \(([\cdot, \cdot]_\mu, \omega_\mu, \check{g}_\mu)(t)\) gives us a well-defined flow of Riemannian metrics \(\tilde{g}_{\mu(t)}\) on a flow of trivial principal bundles \(\tilde{M}_{\mu(t)} \to \RR\). 
    As \(\delta_\mu(L) = 0\), the bracket on the adjoint bundle is parallel with respect to the metric connection, that is, in a parallel frame, the bracket is constant. 
    Because of this, we can write the bracket flow as 
    \[
    (\tilde{M}_\infty, \tilde{g}_{\mu(t)}, \tilde{p}) \cong (\RR \times \mathsf{N}_{\mu(t)}, ds^2 +h_{\mu(t)}(s), (0, \mathsf{e})),
    \]
     where \(h_{\mu(t)}(s)\) is a flow of one-parameter families of right-invariant metrics on the flow of Lie groups \(\mathsf{N}_{\mu(t)}\). 
     Let \(t_k \to \infty\) be some given sequence and take \(\mu_k := \mu(t_k)\).
    By \autoref{prop:bracket-bound}, this sequence lies within a bounded set in a finite-dimensional space.
     After passing to a subsequence, \(\mu_k\) converges to a      nilpotent bracket \(\mu_\infty\). 
    The Lie exponential maps \(\exp_{\mu_k}:(\mathfrak{n}_\infty, \mu_k) \to \mathsf{N}_{\mu_k}\)  give us the Cheeger--Gromov  diffeomorphisms:
    \[
    \Phi_k: \RR \times \mathsf{N}_{\mu_\infty} \to \RR \times \mathsf{N}_{\mu_k}, 
    \qquad
    \Phi_k(s, \exp_{\mu_\infty}(X)) = (s, \exp_{\mu_k}(X)). 
    \] 
    Indeed, \(\Phi_k(0, \mathsf{e}) = (0, \mathsf{e})\) and
    \[
    \Phi_k^*(ds^2 +h_{\mu_k}(s)) \to ds^2 +h_{\mu_\infty}(s)
    \]
    smoothly on compact sets, 
    where \(h_{\mu_\infty}(s)\) is the right-invariant metric on \(\mathsf{N}_{\mu_\infty}\) which agrees with \(h_{\mu_k}(s)\) along the identities. 
    By \autoref{lemma:Lauret-monotone}, the fibres \((\mathsf{N}_{\mu_\infty}, h_{\mu_\infty}(s))\) are nilsolitons with soliton constant \(-\frac{1}{2}\). 
    Call the soliton derivations \(D_\infty(s) \in {\rm Der}(\mathfrak{n}_\infty, \mu_\infty)\). 
    We retain the properties
    \[
    \delta_{\mu_\infty}( L_\infty )=0, 
    \qquad |L_\infty|^2 = \frac{1}{2}, 
    \qquad {\hat\nabla}^*L_\infty = 0, \qquad \nabla L_\infty = 0. 
    \]
     From the first property, \(L_\infty \in {\rm Der}(\mathfrak{n}_\infty, \mu_\infty)\).
     In the horizontal trivialisation \(\RR \times \mathsf{N}_{\mu_\infty}\), the shape operator \(L_\infty\) is constant, hence 
    \[
    h_{\mu_\infty}(s) = (e^{sL_\infty})^* h_{\mu_\infty}(0).
    \]  
    Indeed, \(g_{\mu_\infty}=ds^2 +h_{\mu_\infty}(s)\) is a homogeneous metric that is invariant under right multiplication by the solvable group 
    \[
    \mathsf{S}_\infty = \RR \ltimes_{L_\infty} \mathsf{N}_{\mu_\infty}, 
    \qquad [a, U]  = L_\infty U.
    \]
    Here the multiplication is \((s, \mathsf{n})(r, \mathsf{m}) =(s+r, \phi_r(\mathsf{n})\mathsf{m})\), where \(d\phi_r = e^{-rL_\infty}\). 
    It remains to show that \(g_{\mu_\infty}\) is a solvsoliton. 
    Consider O'Neill's formulas for the projection \(\mathsf{S}_\infty \to \RR\) (\autoref{prop:Ricci-curv-dim1}):
    \begin{align*}
        \Ric_{g_{\mu_\infty}}^{\rm ad}(U) &=  -\frac{1}{2} U + D_\infty(s)U -(\tr L_\infty)L_\infty U,
        \\
        \Ric_{g_{\mu_\infty}}^{\mathrm{ad}, TB}(X)&= 0,
        \\
        \Ric_{g_{\mu_\infty}}^{TB}(X) &= -\frac{1}{2}X.
    \end{align*}
    The nilsoliton derivation is parallel since, by \autoref{prop:1st-2nd-variation-moment-map},
    \[
    h_{\mu_\infty}(\nabla D_\infty, \alpha) = h_{\mu_\infty}(\nabla\hat{\Ric}_{\mu_\infty}, \alpha) = \frac{1}{2} h_{\mu_\infty}(\delta_{\mu_\infty}(\alpha), \delta_{\mu_\infty}(L_\infty)) = 0
    \]
    for all \(\alpha \in {\rm End}(\mathfrak{n}_\infty)\). 
    As \(L_\infty\) is a symmetric derivation, it commutes with \(\hat\Ric_{\mu_\infty}\), and thus \([D_\infty, L_\infty] = 0\) \cite[Lem. 2.2]{heber1998noncompact}.
    Indeed, \(g_{\mu_\infty}\) is a solvsoliton with constant \(-\frac{1}{2}\) and solvsoliton derivation equal to
    \[
    D_\infty -(\tr L_\infty)L_\infty
    \]    
    on \(\mathfrak{n}_\infty\) and zero on \(\RR a\). 
\end{proof}

\subsection{Base collapse}

\begin{coro}[\autoref{coro:D}]
\label{coro:D-precise}
    Let \(\mathsf{N}\) be a simply-connected nilpotent Lie group. 
    Let \(M \to S^1\) be a twisted principal \(\mathsf{N}\)-bundle with unipotent monodromy. 
    If \(g(t)\) is a normalised Ricci flow of invariant metrics on \(M\), then
    \[
    \lim_{t \to \infty}\diam(S^1, \check{g}(t)) =  0,
    \]
    where \(\check{g}(t)\) is the quotient metric.
\end{coro}

\begin{proof}
    The length of the base is \(\ell(t) = 2\diam(S^1, \check{g}(t))\).
    By \autoref{lemma:four-properties}, \(\ell(t) \to \ell_\infty \geq 0\). 
    Suppose for contradiction that \(\ell_\infty> 0\). 
    Then, by \autoref{lemma:four-properties}, there exists a constant \(d>0\) such that 
    \[
    0< d\leq \diam(S^1, \check{g}(t))
    \]
    for all \(t \geq 0\). 
    By \autoref{thm:inv-RF-compactness}, the time-translated normalised Ricci flows \(g(t+ t_k)\) subconverge to an invariant normalised Ricci flow on a principal \(\mathsf{N}_\infty\)-bundle over \(\RR\), 
    \[
    (\tilde{M}_\infty, \tilde{g}_\infty(t), \tilde{p}_\infty), 
    \]
    where \(\mathsf{N}_\infty\) is a simply connected nilpotent Lie group. 
    From the diameter bound, \((\tilde{M}_\infty, \tilde{g}_\infty(t))\) descends to an invariant normalised Ricci flow on a twisted principal \(\mathsf{N}_\infty\)-bundle \(M_\infty \to S^1\).
    By \autoref{lemma:four-properties}, 
    \((\tilde{M}_\infty, \tilde{g}_\infty(t))\) has the following properties:
    \begin{equation}\label{equ:corollary-pf-properties}
        \delta_{[\cdot, \cdot]}( L_\infty)=0, 
    \qquad |L_\infty|^2 = \frac{1}{2}, 
    \qquad {\hat\nabla}^*L_\infty = 0, \qquad \nabla L_\infty = 0. 
    \end{equation}
    From the last property, \(L_\infty\) is \(\nabla^{\rm ad}\)-parallel. 
    In a horizontal trivialisation, \(L_\infty\) is constant, so if \(h_\infty: \RR \to {\rm Sym}^2_+(\mathfrak{n}_\infty^*)\) is the vector bundle metric on the adjoint bundle, then
    \[
    h_\infty(s) = (e^{sL_\infty})^*h_\infty(0).
    \]
    If \( P = d\varphi_\mathsf{e} \circ \operatorname{Ad}(\mathsf{n})\)
    is the adjoint holonomy of \(\nabla^{\rm ad}\) at time \(t = 0\), then 
    \begin{equation}\label{equ:adjoint-holonomy}
        P^*P = e^{2\ell_\infty L_\infty}, 
    \end{equation}
    where \(P^*\) is the \(h_\infty(0)\)-adjoint of \(P\). 
    Since \(L_\infty\) is \(\nabla^{\rm ad}\)-parallel, it commutes with \(P\). 
    It follows by \eqref{equ:adjoint-holonomy} that \(P\) commutes with \(P^*P\), hence 
    \[
    P^*P  = PP^*, 
    \]
    that is, \(P\) is normal. 
    \newline
    \indent
    %
    Let \(P_k \in {\rm Aut}(\mathfrak{n})\) be the adjoint holonomy of \(\nabla^{\rm ad}_k\) at time \(t = 0\). 
    Write \(P_k = d(\varphi_k)_\mathsf{e} \circ \operatorname{Ad}(\mathsf{n}_k)\), where \(\varphi_k\) is the monodromy and \(\mathsf{n}_k \in \mathsf{N}\). 
    Let \(\{\mathfrak{n}^j\}\) denote the lower central series of \(\mathfrak{n}\). 
    As the exponential map is surjective, there exists \(X_k \in \mathfrak{n}\) such that \(\operatorname{Ad}(\mathsf{n}_k) = \exp (\operatorname{ad}_{X_k})\). 
    Thus, 
    \[
    (\operatorname{Id}-\operatorname{Ad}(\mathsf{n}_k))(\mathfrak{n}^j) \subseteq \mathfrak{n}^{j+1},
    \]
    and as automorphisms preserve the lower central series, \(P_k = d(\varphi_k)_\mathsf{e} \circ \operatorname{Ad}(\mathsf{n}_k) = d(\varphi_k)_\mathsf{e}\) on the quotients \(\mathfrak{n}^{j}/\mathfrak{n}^{j+1}\). 
    By assumption, \(d(\varphi_k)_\mathsf{e}\) is unipotent, so  \((\operatorname{Id} - d(\varphi_k)_\mathsf{e})^n = 0\) where \(n = \dim \mathfrak{n}\). 
    Hence,
    \[
    (\operatorname{Id}-P_k)^n(\mathfrak{n}^j)=(\operatorname{Id}-d(\varphi_k)_\mathsf{e}\circ\operatorname{Ad}(\mathsf{n}_k))^n(\mathfrak{n}^j) \subseteq \mathfrak{n}^{j+1},
    \]
    so if \(\mathfrak{n}\) is \(m\)-step nilpotent, then 
    \[
      (\operatorname{Id}-P_k)^{nm}(\mathfrak{n})= (\operatorname{Id}-d(\varphi_k)_\mathsf{e}\circ\operatorname{Ad}(\mathsf{n}_k))^{nm}(\mathfrak{n}) \subseteq \mathfrak{n}^{m+1} = \{0\}.
    \]
    Indeed, \(P_k\) is unipotent, and by \autoref{prop:adjoint-holonomy-compactness}, after passing to a subsequence, there exist linear isometries \(\rho_k : \mathfrak{n}_\infty \to \mathfrak{n}\) such that \(\rho_k^{-1} P_k \rho_k \to P\). 
    It follows that 
    \[
    (\operatorname{Id} - P)^n = \lim_{k \to \infty} \rho_k^{-1} (\operatorname{Id} - P_k)^n \rho_k = 0,
    \]
    and therefore, \(P\) is unipotent. 
    However, a unipotent normal matrix is the identity, so \eqref{equ:adjoint-holonomy} becomes
    \[
    \operatorname{Id} = e^{2\ell_\infty L_\infty}.
    \] 
    This means \(L_\infty = 0\), which contradicts \(|L_\infty|^2 = \frac{1}{2}\) in \eqref{equ:corollary-pf-properties}.
\end{proof}

\bibliographystyle{alpha}
\bibliography{bibl}

\end{document}